\documentclass[12pt]{article}
\usepackage{amsfonts,amsmath,amssymb,amsthm}
\usepackage{enumerate}
\usepackage{xparse}
\usepackage{ifthen}
\usepackage{float}
\usepackage{tikz}
\usetikzlibrary{calc, shadows}
\usepackage[margin=1in]{geometry}

\theoremstyle{definition}
	\newtheorem{defn}{Definition}[section]
\theoremstyle{definition}
	
\theoremstyle{plain}
	\newtheorem{thm}[defn]{Theorem}
	\newtheorem{lem}[defn]{Lemma}
	\newtheorem{cor}[defn]{Corollary}
	\newtheorem{prop}[defn]{Proposition}
\theoremstyle{remark}
	\newtheorem{rem}{Remark}

\newcommand{\ignore}[1]{}
\newcommand{\etal}{\mbox{\textit{et al}.\ }}

\newcommand{\ii}{\mathtt{i}}

\newcommand{\Z}{\mathbb{Z}}
\newcommand{\Q}{\mathbb{Q}}

\newcommand{~}{\sim}

\newcommand{\supp}{\mathtt{supp}}

\DeclareMathOperator{\Sp}{Sp}

\DeclareMathOperator{\Aut}{Aut}

\newcommand{\Damp}{\Lambda}

\newcommand{\emptyg}[1]{\overline{K}_{#1}}
\newcommand{\norm}[1]{\left\lVert #1 \right\rVert}
\newcommand{\allone}[1]{{j}_{#1}}

\newcommand{\twovector}[2]{\left[\begin{array}{c} #1 \\ #2 \end{array}\right]}

\newcommand{\corona}{\circ}

\NewDocumentCommand\eigenA{mg}{%
    \ensuremath{E_{#1}{\IfNoValueTF{#2}{}{(#2)}}}%
}
\NewDocumentCommand\eigenL{mg}{%
    \ensuremath{F_{#1}{\IfNoValueTF{#2}{}{(#2)}}}%
}

\ifthenelse{\isundefined{\usedirac}}
{ 
        \newcommand{\ket}[1]{\mathbf{#1}}
        \newcommand{\bra}[1]{\mathbf{#1}^{\dagger}}
        \newcommand{\tbra}[1]{\mathbf{#1}^{T}}
        \newcommand{\eket}[1]{{\mathbf{e}}_{#1}}        
        \newcommand{\ebra}[1]{\eket{#1}^{T}}            
        \newcommand{\outprod}[2]{\ket{#1}\bra{#2}}                              

}
{ 
        \newcommand{\bra}[1]{\left\langle #1 \right|}
        \newcommand{\ket}[1]{\left| #1 \right\rangle}
        \newcommand{\ebra}[1]{\bra{#1}}
        \newcommand{\eket}[1]{\ket{#1}}

        \newcommand{\outprod}[2]{\ket{#1}\! \bra{#2}}

}

\newcommand{\topline}[1]{}
\newcommand{\botline}[1]{}

\newif\ifnotesw\noteswtrue
   {\ifnotesw\marginpar[\hfill\(\top\)]{\(\top\)}\fi}%
   {\ifnotesw\marginpar[\hfill\(\bot\)]{\(\bot\)}\fi}

\newcommand{\mnote}[1]%
    {\ifnotesw\marginpar%
        [{\scriptsize\begin{minipage}[t]{\marginparwidth}
        \raggedleft#1%
                        \end{minipage}}]%
        {\scriptsize\begin{minipage}[t]{\marginparwidth}
        \raggedright#1%
                        \end{minipage}}%
    \fi}

\begin{document}

\title{Quantum State Transfer in Coronas}
\author{
Ethan Ackelsberg\footnote{Division of Science, Mathematics, and Computing, Bard College at Simon's Rock}
\and
Zachary Brehm\footnote{Department of Mathematics, SUNY Potsdam}
\and
Ada Chan\footnote{Department of Mathematics and Statistics, York University}
\and
Joshua Mundinger\footnote{Department of Mathematics and Statistics, Swarthmore College}
\and
Christino Tamon\footnote{Department of Computer Science, Clarkson University}}
\maketitle 

\begin{abstract}
We study state transfer in quantum walk on graphs relative to the adjacency matrix. 
Our motivation is to understand how the addition of pendant subgraphs affect state transfer.
For two graphs $G$ and $H$, the Frucht-Harary corona product $G \corona H$ is obtained by taking 
$|G|$ copies of the cone $K_{1} + H$ and by connecting the conical vertices according to $G$.
Our work explores conditions under which the corona $G \corona H$ exhibits state transfer.
We also describe new families of graphs with state transfer based on the corona product.
Some of these constructions provide a generalization of related known results.
\end{abstract}


\section{Introduction}

Quantum walk is a natural generalization of classical random walk on graphs.
It has received strong interest due to its important applications in quantum information and computation.
Farhi and Gutmann \cite{fg98} introduced quantum walk algorithms for solving search problems on graphs. 
In their framework, given a graph $G$ with adjacency matrix $A$, the time-varying unitary matrix 
$U(t) := e^{-\ii tA}$ defines a continuous-time quantum walk on $G$.
Subsequently Childs \etal \cite{ccdfgs03} showed that these algorithms may provide exponential speedup 
over classical probabilistic algorithms.
The work by Farhi \etal \cite{fgg08} described an intriguing continuous-time quantum walk algorithm 
with nontrivial speedup for a concrete problem called Boolean formula evaluation.

In quantum information, Bose \cite{b03} studied the problem of information transmission in quantum spin chains.
Christandl \etal \cite{cdel04,cddekl05} showed that this problem may be reduced to the following phenomenon in
quantum walk. Given two vertices $u$ and $v$ in $G$, 
we say perfect state transfer occurs between $u$ and $v$ in $G$ at time $\tau$ if the $(u,v)$-entry 
of $U(\tau)$ has unit magnitude. 
It was quickly apparent that perfect state transfer is an exotic phenomenon.
Godsil \cite{g12} proved that for every positive integer $k$ 
there are only finitely many graphs of maximum degree $k$ which have perfect state transfer. 
Therefore, the following relaxation of this notion is often more useful to consider.
The state transfer between $u$ and $v$ is called 
``pretty good'' (see Godsil \cite{g11-dm}) or
``almost perfect'' (see Vinet and Zhedanov \cite{vz12}) 
if the $(u,v)$-entry of the unitary matrix $U(t)$ can be made arbitrarily close to one.

Christandl \etal \cite{cdel04, cddekl05} observed that the path $P_{n}$ on $n$ vertices
has antipodal perfect state transfer if and only if $n=2,3$. 
In a striking result, Godsil \etal \cite{gkss12} proved that $P_{n}$ has antipodal pretty good state transfer 
if and only if $n+1$ is a prime, twice a prime, or a power of two.
This provides the first family of graphs with pretty good state transfer which
correspond to the quantum spin chains originally studied by Bose.

Shortly after,
Fan and Godsil \cite{fg13} studied a family of graphs obtained by taking two cones $K_{1} + \emptyg{m}$ and
then connecting the two conical vertices. They showed that these graphs, which are called double stars,
have no perfect state transfer, but have pretty good state transfer between the two conical vertices 
if and only if $4m+1$ is a perfect square.
These graphs provide the second family of graphs known to have pretty good state transfer.

In this work, we provide new families of graphs with pretty good state transfer.
Our constructions are based on a natural generalization of Fan and Godsil's double stars.
The corona product of an $n$-vertex graph $G$ with another graph $H$, typically denoted $G \corona H$, 
is obtained by taking $n$ copies of the cone $K_{1} + H$ and by connecting the conical vertices according to $G$.
In a corona product $G \corona H$, we sometimes call $G$ the {\em base} graph and $H$ the {\em pendant} graph.
This graph product was introduced by Frucht and Harary \cite{fh70} in their study of 
automorphism groups of graphs which are obtained by wreath products.

We first observe that perfect state transfer on corona products is extremely rare.
This is mainly due to the specific forms of the corona eigenvalues 
(which unsurprisingly resemble the eigenvalues of cones)
coupled with the fact that periodicity is a necessary condition for perfect state transfer.
Our negative results apply to corona families $G \corona H$ when $H$ is either the empty or the complete graph, 
under suitable conditions on $G$. 
In a companion work \cite{abcmt-laplacian}, we observed an optimal negative result which holds for all $H$ 
but in a Laplacian setting.

Given that perfect state transfer is rare, our subsequent results mainly focus on pretty good state transfer.
We prove that the family of graphs $K_{2} \corona K_{m}$,
which are called barbell graphs (see Ghosh \etal \cite{gbs08}), 
admit pretty good state transfer for all $m$. 
Here, state transfer occurs between the two vertices of $K_{2}$.
This is in contrast to the double stars $K_{2} \corona \emptyg{m}$ where pretty good state transfer
requires number-theoretic conditions on $m$.

We observe something curious for corona products when the base graph is complete.
It is known that the complete graphs $K_{n}$, $n\geq 3$, are periodic but have no perfect state transfer, 
and hence have no pretty good state transfer.
We observe that although $K_{n} \corona \emptyg{2}$ has no pretty good state transfer,
its minor variant $(K_{n} \Box K_{2}) \corona \emptyg{2}$ has pretty good state transfer 
for all but two values of $n$. Here, the state transfer occurs between the two vertices of $\emptyg{2}$.
An immediate corollary is that $(K_{n} + I) \corona \emptyg{2}$ has pretty good state transfer, 
where $K_{n} + I$ denotes the graph obtained by adding self-loops to each vertex of $K_{n}$.
This provides another example where self-loops are useful for state transfer 
(see Casaccino \etal \cite{clms09}).

Our other results involve graphs of the form $G \corona K_{1}$ which are called thorny graphs (see Gutman \cite{g98}).
We show that if $G$ is a graph with perfect state transfer at time $\pi/g$, for some positive integer $g$,
then $G \corona K_{1}$ has pretty good state transfer, provided that the adjacency matrix of $G$ is nonsingular.
On the other hand, if the adjacency matrix of $G$ is singular, we derive the same result 
if $G$ has perfect state transfer at time $\pi/2$.
Taken together, this proves that the thorny cube $Q_{d} \corona K_{1}$ has pretty good state transfer for all $d$.
This confirms some of the numerical observations of Makmal \etal \cite{mzmtb14} on embedded hypercubes
albeit for the continuous-time setting (since their results were stated for discrete-time quantum walks).

It turns out that perfect state transfer is not necessary for a thorny graph to have state transfer.
Coutinho \etal \cite{cggv15} proved that the cocktail party graph $\overline{nK_{2}}$ has perfect state transfer 
if and only if $n$ is even.
We show that the thorny graph $\overline{nK_{2}} \corona K_{1}$ has pretty good state transfer for all $n$.
This observation holds in the Laplacian setting as well (see \cite{abcmt-laplacian}).

Returning to the double stars,
we also observe that $C_{2} \corona \emptyg{m}$ has {\em perfect} state transfer 
whenever $m+1$ is an even square.
Here, $C_{2}$ is the digon (which is a multigraph on two vertices 
connected by two parallel edges). This shows that certain double stars are one additional edge away 
from having perfect state transfer. This is akin to using weighted edges on paths for perfect state transfer 
(see Christandl \etal \cite{cdel04}). Here, we merely use an integer weight in a non-unimodal
weighting scheme on a path.


\section{Preliminaries}

Given an undirected graph $G = (V,E)$ on $n$ vertices, the adjacency matrix of $G$ is an $n \times n$ 
symmetric matrix $A(G)$ where, for all vertices $u$ and $v$, the $(u,v)$ entry of $A(G)$ is $1$, 
if $(u,v) \in E$, and $0$, otherwise.
The spectrum of $G$, which we denote as $\Sp(G)$, is the set of distinct eigenvalues of $A(G)$.
We use $\rho(G)$ to denote the spectral radius of $G$.
By the spectral theorem, we may write
\begin{equation}
A(G) = \sum_{\lambda \in \Sp(G)} \lambda \eigenA{\lambda}{G}
\end{equation}
where $\eigenA{\lambda}{G}$ is the eigenprojector (orthogonal projector onto an eigenspace)
corresponding to eigenvalue $\lambda$.

The {\em eigenvalue support} of a vertex $u$ in $G$, denoted $\supp_{G}(u)$,
is the set of eigenvalues $\lambda$ of $G$ for which $\eigenA{\lambda}{G}\eket{u} \neq 0$. 
Two vertices $u$ and $v$ of $G$ are called {\em strongly cospectral} if 
$\eigenA{\lambda}{G}\eket{u} = \pm\eigenA{\lambda}{G}\eket{v}$ 
for each eigenvalue $\lambda$ of $G$.

A continuous-time quantum walk on $G$ is defined by the time-varying unitary matrix 
$U(t) = e^{-\ii tA(G)}$. 
We say such a quantum walk has {\em perfect state transfer} 
between vertices $u$ and $v$ at time $\tau$ if 
\begin{equation}
\ebra{v}U(\tau)\eket{u} = \gamma,
\end{equation} 
for some complex unimodular $\gamma$. 
We call $\gamma$ the {\em phase} of the perfect state transfer.
Note that $U(t)$ is symmetric for all $t$.
In the special case where perfect state transfer occurs with $u = v$,
we say the quantum walk is {\em periodic at vertex $u$}. Moreover, if $U(\tau)$ is a scalar
multiple of the identity matrix, then the quantum walk is {\em periodic}.

Given that perfect state transfer is rare,
we consider the following relaxation of this phenomenon proposed by Godsil. 
The quantum walk on $G$ has {\em pretty good state transfer} between $u$ and $v$ if
for all $\epsilon > 0$ there is a time $\tau$ so that
\begin{equation}
\norm{e^{-\ii\tau A(G)}\eket{u} - \gamma\eket{v}} < \epsilon,
\end{equation}
for some complex unimodular $\gamma$.

In what follows, we state some useful facts about state transfer on graphs.
The following formulation by Coutinho (which summarized the most relevant facts) will be useful
for our purposes.

\begin{thm}[Coutinho \cite{c14}, Theorem 2.4.4] \label{thm: coutinho-pst-conditions}
Let $G$ be a graph and let $u,v$ be two vertices of $G$.  
Then there is perfect state transfer between $u$ and $v$
at time $t$ with phase $\gamma$ if and only if all of the following conditions hold.
\begin{enumerate}[i)]
\item Vertices $u$ and $v$ are strongly cospectral.

\item There are integers $a,\Delta$ where $\Delta$ is square-free so that 
	for each eigenvalue $\lambda$ in $\supp_{G}(u)$:
	\begin{enumerate}[(a)]
	\item $\lambda = \frac{1}{2}(a + b_{\lambda}\sqrt{\Delta})$,
		for some integer $b_{\lambda}$.
	\item $\ebra{u}\eigenA{\lambda}{G}\eket{v}$ is positive if and only if
		$(\rho(G) - \lambda)/g\sqrt{\Delta}$ is even,
		where 
		\begin{equation}
			g := \gcd\left( \left\{ \frac{\rho(G) - \lambda}{\sqrt{\Delta}} : 
				\lambda \in \supp_{G}(u) \right\} \right).
		\end{equation}
	\end{enumerate}
\end{enumerate}
Moreover, if the above conditions hold, then the following also hold.
\begin{enumerate}[i)]
\item There is a minimum time of perfect state transfer between $u$ and $v$ given by
	\begin{equation}
	t_0 := \frac{\pi}{g\sqrt{\Delta}}.
	\end{equation}

\item The time of perfect state tranfer $t$ is an odd multiple of $t_0$.

\item The phase of perfect state transfer is given by $\gamma = e^{-\ii t\rho(G)}$.
\end{enumerate}
\end{thm}

We state some strong properties proved by Godsil which relate perfect state transfer with periodicity
in a fundamental way.

\begin{lem}[Godsil \cite{g11}] \label{lem: periodic-necessary-condition}
If $G$ has perfect state transfer between vertices $u$ and $v$ at time $t$,
then $G$ is periodic at $u$ at time $2t$.
\end{lem}

\begin{thm}[Godsil \cite{g12}] \label{thm: godsilperiodic}
A graph $G$ is periodic at vertex $u$ if and only if either:
\begin{enumerate}[i)]
\item all eigenvalues in $\supp_{G}(u)$ are integers; or
\item there is a square-free integer $\Delta$ and an integer $a$ so that
	each eigenvalue $\lambda$ in $\supp_{G}(u)$ is of the form 
	$\lambda = \frac{1}{2}(a + b_{\lambda}\sqrt{\Delta})$, 
	for some integer $b_{\lambda}$.
\end{enumerate}
\end{thm}

We also state some necessary conditions for perfect state transfer and  pretty good state transfer in terms of the automorphisms of $G$.

\begin{thm}[Godsil \cite{g12}] \label{thm: godsil-automorphism}
Let $G$ be a graph with perfect state transfer between vertices $u$ and $v$.
Then, for each automorphism $\tau \in \Aut(G)$, $\tau(u) = u$ if and only if $\tau(v) = v$. 
\end{thm}

\begin{thm}[Godsil \cite{g15}, Lemmas 7.4.1 and 9.1.4]  \label{thm:pgst-necessary-conditions}
Suppose $G$ has pretty good state transfer between vertices $u$ and $v$.
Then $u$ and $v$ are strongly cospectral, and each automorphism fixing $u$ must fix $v$.
\end{thm}

\bigskip

In our analysis, we will need some tools from number theory.
For example, we will need the following form of Kronecker's approximation theorem.

\begin{thm}[Kronecker's Theorem: Hardy and Wright \cite{hw00}, Theorem 442] \label{thm: kronecker} 
Let $1,\lambda_{1},\ldots,\lambda_{m}$ be linearly independent over $\Q$.
Let $\alpha_{1},\ldots,\alpha_{m}$ be arbitrary real numbers,
and let $N,\epsilon$ be positive real numbers.
Then there are integers $\ell > N$ and $q_1,\ldots, q_m$
so that
\begin{equation}\label{eq: kroneckers-theorem}
|\ell\lambda_{k} - q_{k} - \alpha_{k}| < \epsilon,
\end{equation}
for each $k=1,\ldots,m$.
\end{thm}

For brevity, whenever we have an inequality of the form $|\alpha - \beta| < \epsilon$, for arbitrarily small $\epsilon$,
we will write instead $\alpha \approx \beta$ and omit the explicit dependence on
$\epsilon$. For example, \eqref{eq: kroneckers-theorem} will be represented as 
$\ell\lambda_k - q_k \approx \alpha_k$.

\medskip

In our applications of Kronecker's Theorem, we will use the following lemma to identify 
sets of numbers which are linearly independent over the rationals.

\begin{lem}[Newman and Flanders \cite{f60}] \label{lem: integersquareroot}
Let $a_{1},\ldots,a_{n}$ be positive integers, coprime in pairs, no one of which is a perfect square.
Then the $2^{n}$ algebraic integers
\begin{equation}
	\sqrt{a_{1}^{e_{1}} \ldots a_{n}^{e_{n}}},
	\hspace{1.5in}
	e_{j} = 0,1,
\end{equation}
are linearly independent over the field $\Q$ of rationals.
\end{lem}

\paragraph{Notation}
We describe some notation used throughout the rest of the paper.

The all-one and all-zero vectors of dimension $n$ are denoted $\ket{\allone}_{n}$, $\ket{0}_{n}$, respectively.
The $m \times n$ all-one matrix is denoted $J_{m,n}$ or simply $J_{n}$ if $m = n$.
The identity matrix of size $n$ is denoted $I_{n}$.

For standard families of graphs, we use $K_{n}$ for the complete graph on $n$ vertices,
$P_{n}$ for the path on $n$ vertices, and $Q_{n}$ for the $n$-dimensional cube. 
In what follows, let $G$ and $H$ be graphs.
The complement of $G$ is denoted $\overline{G}$.
The disjoint union of $G$ and $H$ is written as $G \cup H$, while the disjoint union
of $n$ copies $G$ is denoted $nG$. 
The Cartesian product of $G$ and $H$ is denoted $G \Box H$.
If $G,H$ have $n,m$ vertices, respectively, the adjacency matrix of $G \Box H$ is given by
$A(G) \otimes I_{m} + I_{n} \otimes A(H)$.
The join $G + H$ of $G$ and $H$ is the graph $\overline{\overline{G} \cup \overline{H}}$.
Further background on algebraic graph theory may be found in Godsil and Royle \cite{gr01}.


\section{Corona of Graphs}

We define the Frucht-Harary corona product of two graphs (see \cite{fh70}).
Let $G$ be a graph on the vertex set $\{v_1,\ldots, v_n\}$ and $H$ be a graph
on the vertex set $\{1,\ldots,m\}$. The latter is chosen for notational convenience.
The \textit{corona} $G \corona H$ is formed by taking the disjoint union of $G$ 
and $n$ copies of $H$ and then adding an edge from each vertex of the $j$th copy of $H$ to
the vertex $v_j$ in $G$. 
Formally, the corona $G \corona H$ has the vertex set
\begin{equation}
	V(G \corona H) = V(G) \times \left( \{0\} \cup V(H) \right), 
\end{equation}
and the adjacency relation
\begin{equation}
	((v, w), (v', w')) \in E(G \corona H)
	\iff
	 \begin{cases}
	 	\text{$w = w' = 0$ and $(v,v') \in E(G)$} & \text{or} \\
		\text{$v = v'$ and $(w,w') \in E(H)$} & \text{or} \\
		\text{$v = v'$ and exactly one of $w$ and $w'$ is $0$.}
	\end{cases}
\end{equation}

\begin{figure}[H]
	\centering
	\topline{.15in}
	\begin{tikzpicture}[scale=1.25]
	\filldraw[fill=white, circular drop shadow] (0,0) node{$G$} circle[x radius = 1.7, y radius = 0.7];
	\foreach \x in {1,2,4,5} {
		\foreach \y in {0,...,5}
			\draw (-1.5+0.5*\x, -0.1*\x*\x + 0.6*\x - 0.5) -- ++(-2+0.5*\x+0.2*\y, 1);
		\filldraw[fill=white, circular drop shadow] (-1.5+0.5*\x, -0.1*\x*\x + 0.6*\x - 0.5) circle[x radius = 0.03in, y radius = 0.02in] node{};
	}
	\filldraw[fill=white, circular drop shadow] (-2, 1) circle[x radius=0.5, y radius=0.3] node{$H$};
	\filldraw[fill=white, circular drop shadow] (-1, 1.3) circle[x radius=0.5, y radius=0.3] node{$H$};
	\draw[thick, loosely dotted] (-0.3,1.3) -- ++(0.6,0);
	\filldraw[fill=white, circular drop shadow] (1, 1.3) circle[x radius=0.5, y radius=0.3] node{$H$};
	\filldraw[fill=white, circular drop shadow] (2, 1) circle[x radius=0.5, y radius=0.3] node{$H$};
    \end{tikzpicture} 
    \caption{A corona $G \corona H$.}
	\botline{.15in}
\end{figure}
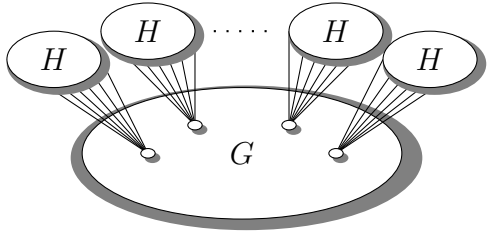

The adjacency matrix of $G \corona H$ where $H$ has $m$ vertices is given by
\begin{equation} \label{eq: corona-adjacency-matrix}
	A(G \corona H) = A(G) \otimes \ebra{0}\eket{0} + I_n \otimes
		\left[\begin{array}{cc}
			0 & \tbra{j}_m \\
			\ket{j}_m & A(H)
		\end{array}\right].
\end{equation}

The spectrum of $G \corona H$ when $H$ is a $k$-regular graph is known and we state
this in the following. 

\begin{thm}[Barik \etal \cite{bps07}] \label{thm: regular-corona-spectrum}
Suppose that $G$ is a graph and $H$ is a $k$-regular graph on $m$ vertices. 
Suppose $G$ has eigenvalues $\lambda_0 > \ldots > \lambda_p$ with multiplicities $r_0, \ldots, r_p$, 
and $H$ has eigenvalues $k = \mu_0 > \mu _1 > \ldots > \mu_q$ with multiplicities $s_0, \ldots, s_q$. 
Then $G \corona H$ has the following spectrum.
\begin{enumerate}[i)]
	\item $k$ is an eigenvalue with multiplicity $n(s_0 - 1)$.
	\item $\mu_j$ is an eigenvalue with multiplicity $ns_j$ for $j = 1, \ldots, q$.
	\item $\frac{1}{2}(\lambda_j + k \pm \sqrt{(\lambda_j - k)^2 + 4m})$ are eigenvalues
		each with multiplicity $r_j$ for $j = 0, \ldots, p$.
\end{enumerate}
\end{thm}

Theorem \ref{thm: regular-corona-spectrum} gives the spectrum of $G \corona H$ when $H$ is $k$-regular. 
We now provide the corresponding eigenprojectors which will be crucial in our subsequent analyses.

\begin{prop} \label{prop: adj-spec-decomp-corona}
Let $G$ be a graph on $n$ vertices and let $H$ be a $k$-regular graph on $m$ vertices.
Then the eigenprojectors of $G \corona H$ are given by the following.
\begin{enumerate}[i)]
\item 
For each eigenvalue $\mu$ of $H$, let
\begin{equation}
	\eigenA{\mu}
	:=  I_n \otimes
	\left[\begin{array}{cc}
		0 & \tbra{0}_m \\ \ket{0}_m & \eigenA{\mu}{H} - \delta_{\mu,k}\frac{1}{m}J_{m} 
	\end{array}\right].
\end{equation}

\item 
For each eigenvalue $\lambda$ of $G$, define a pair of eigenvalues
\begin{equation}
	\lambda_{\pm} := \frac{\lambda + k \pm \sqrt{(\lambda - k)^{2} + 4m}}{2},
\end{equation}
and let
\begin{equation} \label{eq: adj-lambda-pm-projector}
	\eigenA{\lambda_\pm}
	:= \eigenA{\lambda}{G} \otimes \frac{1}{(\lambda_\pm-k)^2 + m} 		
	\left[ \begin{array}{cc} 
   			(\lambda_{\pm}-k)^2 & (\lambda_\pm - k)\tbra{\allone}_{m} \\
   			(\lambda_\pm - k)\ket{\allone}_{m} & J_{m}
		\end{array} \right].
\end{equation}
\end{enumerate}
Then the spectral decomposition of the corona $G \corona H$ is given by
\begin{equation}
	A(G \corona H) = 
	\sum_{\lambda \in \Sp(G)} \sum_\pm \lambda_{\pm} \eigenA{\lambda_\pm}
	 + \sum_{\mu \in \Sp(H)} \mu \eigenA{\mu}.
\end{equation}

\begin{proof}
Recall that the adjacency matrix of $G \corona H$ is given in \eqref{eq: corona-adjacency-matrix}.

For an eigenvalue $\mu$ of $H$,
let $B_{\mu}$ be an orthonormal basis of the $\mu$-eigenspace that is orthogonal to $\ket{\allone}_{m}$. 
If $H$ is a connected $k$-regular graph, then $B_{k}$ is empty.
Suppose $\mu$ has an eigenvalue multiplicity of $r$.
Then the cardinality of $B_{\mu}$ is $r$ if $\mu < k$, and is $r-1$ if $\mu=k$.
Moreover, $\mu$ is also an eigenvalue of $G \corona H$ since for each $\ket{x} \in B_{\mu}$
and for every $\ket{y}$ in $\mathbb C^n$, we have
\begin{equation}
	A(G \corona H) \ket{y} \otimes \twovector{0}{\ket{x}} = \mu \ket{y} \otimes \twovector{0}{\ket{x}}.
\end{equation}
An orthonormal basis for the $\mu$-eigenspace of $G \corona H$
is given by
\begin{equation} \label{eq: basisformuspace}
	\left\{\ket{v} \otimes \twovector{0}{\ket{x}} : v \in V(G), \ket{x} \in B_{\mu} \right\}.
\end{equation}
The cardinality of this set is $nr$ if $\mu < k$, and is $n(r-1)$ if $\mu = k$.
These match the multiplicities given in Theorem \ref{thm: regular-corona-spectrum}.
Note that the eigenprojector of $\mu$ in $H$ is given by
\begin{equation}
\eigenA{\mu}{H} = 
	\left\{\begin{array}{ll}
	\sum_{\ket{x} \in B_{\mu}} \outprod{x}{x}	& \mbox{ if $\mu < k$ } \\[1.3ex]
	\frac{1}{m}J_{m} + \sum_{\ket{x} \in B_{\mu}} \outprod{x}{x}	& \mbox{ if $\mu = k$ }
	\end{array}\right.
\end{equation}
Therefore, the eigenprojector of $\mu$ in $G\corona H$ is given by
\begin{equation} \label{eq: coronamuprojector}
\eigenA{\mu}(G \corona H) = 
	I_n \otimes 
		\left[ 
			\begin{array}{cc} 
			0 & \tbra{0}_m \\ 
			\ket{0}_m & \eigenA{\mu}{H} - \delta_{\mu,k}\frac{1}{m}J_{m}
			\end{array}
		\right] .
\end{equation}

Next, for each eigenvalue $\lambda$ of $G$ with a corresponding (normalized) eigenvector $\ket{x}$, 
we have
\begin{equation}
	A(G \corona H) \left( \ket{x} \otimes \twovector{\lambda_{\pm} - k}{\ket{\allone}_{m}}\right)
	=
	\lambda_{\pm} \ket{x} \otimes \twovector{\lambda_{\pm} - k}{\ket{\allone}_{m}}
\end{equation}
provided that
\begin{equation} \label{eq: lambdapm}
	\lambda_{\pm} := \frac{\lambda + k \pm \sqrt{(\lambda + k)^2 + 4m - 4\lambda k}}{2} = \frac{\lambda + k \pm \sqrt{(\lambda - k)^2 + 4m}}{2}.
\end{equation}
After normalizing the above eigenvectors, we see that
the eigenprojectors of $\lambda_{\pm}$ are given by
\begin{equation} \label{eq: coronapmprojector}
\eigenA{\lambda_\pm}(G \corona H) = 
	\eigenA{\lambda}{G} \otimes 
	\frac{1}{(\lambda_\pm-k)^2 + m} 		
	\left[ \begin{array}{cc} 
   			(\lambda_{\pm}-k)^2 & (\lambda_\pm - k)\tbra{\allone}_{m} \\
   			(\lambda_\pm - k)\ket{\allone}_{m} & J_{m}
	\end{array} \right]
\end{equation}
These provide the remaining $2n$ eigenvectors of $G \corona H$.
\end{proof}

\end{prop}

For state transfer in $G \corona H$, we need to analyze the elements of the transition matrix 
$e^{-\ii tA(G \corona H)}$. 
We use Proposition \ref{prop: adj-spec-decomp-corona} to derive a convenient form for our analysis.

\begin{prop} \label{prop: adj-transition-corona}
Let $G$ be a graph and let $H$ be a $k$-regular graph on $m$ vertices.
For vertices $v$ and $v'$ of $G$, we have
\begin{equation} \label{eq: Atransition}
	\ebra{(v,0)} e^{-\ii tA(G \corona H)} \eket{(v',0)}
	= \sum_{\lambda \in \Sp(G)} e^{-\ii t(\lambda + k)/2} 
		\ebra{v}\eigenA{\lambda}{G}\eket{v'}
		\left( \cos\left(t\Damp_\lambda/2\right) 
			- \frac{(\lambda - k)}{\Damp_\lambda} i\sin\left(t\Damp_\lambda/2 \right) \right). 
\end{equation}
where $\Damp_\lambda = \sqrt{(\lambda - k)^2 + 4m}$.

\begin{proof}
Let $\lambda_{\pm} = \frac{1}{2}(\lambda+k \pm \Damp_{\lambda})$.
From Proposition \ref{prop: adj-spec-decomp-corona}, we get that
\begin{multline} \label{eq: Atransition-general}
\ebra{(v,w)} e^{-\ii tA(G \corona H)} \eket{(v',w')} = 
	\sum_{\lambda \in \Sp(G)} 
		e^{-\ii t(\lambda + k)/2}
		\ebra{v}\eigenA{\lambda}{G}\eket{v'} 
		\left(
		\sum_{\pm} e^{\mp \ii t\Damp_{\lambda}/2} 
		\ebra{w} M_{\lambda_{\pm}} \eket{w'}
		\right) \\
	\ + \ 
		\delta_{v,v'} (1-\delta_{w,0}) (1-\delta_{w',0}) 
		\sum_{\mu \in \Sp(H)} e^{-\ii t\mu}
		\left(\ebra{w}\eigenA{\mu}{H}\eket{w'} - \frac{1}{m}\delta_{\mu,k}\right),
\end{multline}
where
\begin{equation}
M_{\lambda_{\pm}} =
	\frac{1}{(\lambda_{\pm}-k)^{2} + m}
	\begin{bmatrix}
    	(\lambda_{\pm}-k)^2 & (\lambda_\pm - k)\tbra{\allone}_{m} \\
    	(\lambda_\pm - k)\ket{\allone}_{m} & J_{m}
	\end{bmatrix}.
\end{equation}
Therefore, 
\begin{equation} \label{eqn:adjqwalk-on-corona}
	\ebra{(v,0)} e^{-\ii tA(G\corona H)} \eket{(v',0)}
	= \sum_{\lambda \in \Sp(G)} 
		e^{-\ii t(\lambda + k)/2}
		\ebra{v} \eigenA{\lambda}{G} \eket{v'}
			\sum_\pm e^{\mp \ii t\Damp_{\lambda}/2} 
				\frac{(\lambda_\pm - k)^2}{(\lambda_\pm - k)^2 + m}.
\end{equation}
But, the inner summation in (\ref{eqn:adjqwalk-on-corona}) simplifies to
\begin{equation}
\sum_{\pm} e^{\mp \ii t\Damp_{\lambda}/2} 
	\frac{(\lambda_{\pm}-k)^{2}}{(\lambda_{\pm}-k)^{2}+m}
	\ = \
	\cos(t\Damp_{\lambda}/2) - \frac{(\lambda - k)}{\Damp_{\lambda}} i\sin(t\Damp_{\lambda}/2),
\end{equation}
since
$\prod_{\pm} ((\lambda_{\pm} - k)^{2} + m) = m\Damp_{\lambda}^{2}$
and
$\prod_{\pm} (\lambda_{\pm} - k) = -m$.
This proves the claim.
\end{proof}
\end{prop}


\section{Perfect State Transfer}

Lemma \ref{lem: periodic-necessary-condition} shows that periodicity is a necessary condition 
for perfect state transfer. We show that there is no perfect state transfer in most coronas
since their vertices are not periodic. 

\subsection{Conditions on Periodicity}

To investigate periodicity in a corona $G \corona H$, 
the following lemma shows that it is sufficient to consider base vertices of the form $(v,0)$, 
where $v$ is a vertex in $G$.

\begin{lem} \label{lem: periodicinbase}
Let $G$ be a graph and let $H$ be a regular graph.
If $(v,w)$ is periodic in $G \corona H$, then $(v,0)$ is periodic in $G \corona H$. 
\begin{proof}
The eigenvalue support of $(v,0)$ is contained within the eigenvalue support of $(v,w)$,
which is evident from the eigenprojectors given in Proposition \ref{prop: adj-spec-decomp-corona}.
\end{proof}
\end{lem}

Theorem \ref{thm: godsilperiodic} provides necessary and sufficient conditions for periodicity. 
In the following, we show that periodicity in a corona places a strong condition on the eigenvalues 
of the base graph.

\begin{lem} \label{lem: coronaperiodic}
Suppose that $G$ is a connected graph on at least two vertices and $H$ is a $k$-regular graph on $m$ vertices. 
Then a vertex $(v,0)$ is periodic in $G \corona H$ 
if and only if 
there exists a positive square-free integer $\Delta$ 
such that for all eigenvalues $\lambda \in \supp_{G}(v)$, 
both $\lambda - k$ and $\sqrt{(\lambda- k)^2 + 4m}$ are integer multiples of $\sqrt{\Delta}$. 
Moreover, if this holds, we have that $\Delta$ divides $2m$.

\begin{proof}
For each eigenvalue $\lambda$ of $G$, let
\begin{equation} \label{eq: eigenvalue-support}
\lambda_{\pm} := \frac{1}{2} \left( \lambda + k \pm \sqrt{(\lambda - k)^2 + 4m} \right).
\end{equation}
By Proposition \ref{prop: adj-spec-decomp-corona}, the eigenvalue support of $(v,0)$ in $G \corona H$ 
is given by
$\{\lambda_{\pm} : \lambda \in \supp_{G}(v)\}$.
If there is a positive square-free integer $\Delta$ so that for each eigenvalue 
$\lambda$ in the support of $v$ in $G$, both $\lambda - k$ and $\sqrt{(\lambda - k)^2 + 4m}$ 
are integer multiples of $\sqrt{\Delta}$, then $(v,0)$ is periodic in $G \corona H$, 
by Theorem \ref{thm: godsilperiodic}. 
This shows sufficiency.
		
Now for necessity. 
Suppose that $(v,0)$ is periodic in $G \corona H$.
By Theorem \ref{thm: godsilperiodic}, 
either all eigenvalues in the support of $(v,0)$ in $G \corona H$ are integers
or
there is an integer $a$ and a square-free integer $\Delta$ so that 
all eigenvalues in the support of $(v,0)$ in $G \corona H$ are of the form
$\lambda = \frac{1}{2}(a + b_{\lambda} \sqrt{\Delta})$,
for some integers $b_{\lambda}$. 
		
If the eigenvalue support of $(v,0)$ is entirely integers,
then for each eigenvalue $\lambda$ in the support of $v$, we conclude that 
$\lambda + k = \lambda_+ + \lambda_-$ 
and $\sqrt{(\lambda - k)^2 + 4m} = \lambda_+ - \lambda_-$ are integers. 
	
Now suppose there is an integer $a$ and a square-free integer $\Delta > 1$ so that
all eigenvalues $\lambda_{\pm}$ in the support of $(v,0)$ are of the form 
\begin{equation} \label{eq: alternate-form}
	\lambda_{\pm} = \frac{1}{2}(a + b_{\pm}\sqrt{\Delta}),
\end{equation}
where $b_{\pm}$ are integers (which depend on $\lambda_{\pm}$). 
Recall from the proof of Proposition \ref{prop: adj-transition-corona} 
that $(\lambda_+-k)(\lambda_- -k)= -m$. Thus,
\begin{equation}
	-m 
	\ = \ 
	\frac{1}{4}((a-2k)^2 + (b_{+} b_{-})\Delta) 
		+ \frac{1}{4}(a-2k)(b_{+} + b_{-})\sqrt{\Delta}.
\end{equation}
Since $\sqrt{\Delta}$ is not an integer, either $a - 2k$ or $b_{+} + b_{-}$ is zero. 
If $b_{+} = - b_{-}$ then we conclude that $a = \lambda_+ + \lambda_- = \lambda + k$.
This implies that $|\supp_G(v)| = 1$ which contradicts that $G$ is connected on at least two vertices.
Otherwise, we have $a=2k$, from which we conclude that 
$\lambda_\pm = k + \frac{1}{2}b_{\pm}\sqrt{\Delta}$. 
From the definition of $\lambda_{\pm}$ in \eqref{eq: eigenvalue-support}, we have 
\begin{subequations}\begin{align}
	(\lambda_+ - k) + (\lambda_- - k) & = \lambda - k, \\
	(\lambda_+ - k) - (\lambda_- - k) & = \sqrt{(\lambda - k)^2 + 4m}.
\end{align}\end{subequations}
Given the alternate form in \eqref{eq: alternate-form}, this also implies that 
$\lambda - k $ and $\sqrt{(\lambda - k)^2 + 4m}$ are half-integer multiples of $\sqrt{\Delta}$. 
Since their squares are rational algebraic integers, their squares must be integers. 
Thus, 
$\lambda - k$ and $\sqrt{(\lambda - k)^2 + 4m}$ are integer multiples of $\sqrt{\Delta}$.

It remains to show that $\Delta$ divides $2m$. The condition that		
$\sqrt{(\lambda-k)^2  +4m}$ is an integer multiple of $\sqrt{\Delta}$ implies that
\begin{equation}
	\sqrt{\frac{(\lambda-k)^2 + 4m}{\Delta}}
\end{equation} is an integer. Furthermore, since $\Delta$ divides $(\lambda-k)^2$,
it must be that $\Delta$ divides $4m$. Since $\Delta$ is square-free,
$\Delta$ divides $2m$.
\end{proof}
\end{lem}

\smallskip
\begin{rem}
The conditions for periodicity at $(v,w)$ in the corona $G \corona H$ imply that
$\lambda \in k + \Z \sqrt{\Delta}$ for all $\lambda$ in the eigenvalue support of $v$ in $G$.
In particular, $v$ must be periodic in $G$.
\end{rem}

\smallskip
\begin{rem}
	The eigenvalues of $G$ need not be integers for $G \corona H$ to have a periodic vertex.
	Consider $P_{3}$ whose spectrum is $\{0, \pm\sqrt{2}\}$. The eigenvalue support of the middle vertex is 
	$\{\pm\sqrt{2}\}$. 
	Thus we may apply Lemma \ref{lem: coronaperiodic} to see that the middle vertex of $P_3$ is periodic in 
	$P_3 \corona \emptyg{m}$ if and only if $\sqrt{2 + 4m}$ is an integral multiple of $\sqrt{2}$, i.e. if and only if 
	$\sqrt{1 + 2m}$ is an integer. Let us take the specific example of $m = 4$. 
	If we let $v$ denote the center vertex of $P_3$, then the eigenvalue support of $(v,0)$ in 
	$P_3 \corona \emptyg{4}$ is $\{\pm\sqrt{2}, \pm 2 \sqrt{2}\}$, so $(v,0)$ is periodic.
\end{rem}

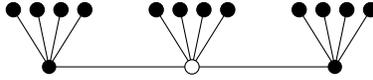
\begin{figure}[h]
	\centering
	\topline{.15in}
    \begin{tikzpicture}[scale=0.95]
	\foreach \x in {-2,0,2} {
		\draw (0,0) -- (\x,0);
		\fill (\x,0) circle[radius = 0.1];
		\foreach \y in {0,...,3}{
			\draw (\x,0) -- +(-1/2 + \y/3, 0.8);
			\filldraw (\x-1/2 + \y/3, 0.8) circle[radius=0.1];
		}	
	}
	\filldraw[fill=white] (0,0) circle[radius=0.1]; 	
	\end{tikzpicture}
    \caption{The white vertex in $P_3 \corona \emptyg{4}$ is periodic.}
	\botline{.15in}
\end{figure}

The heart of our applications of Lemma \ref{lem: coronaperiodic} is the following.
For the conditions of Lemma \ref{lem: coronaperiodic} to be satisfied, it must be that 
$4m/\Delta$ is a difference of squares. 
Since there are only finitely many pairs of squares whose difference is $4m/\Delta$, 
this allows us to restrict possible values of $\lambda - k$ in various ways.
Our first application is to show that the eigenvalues in the support of a periodic vertex 
can not be too close together.

\begin{thm} \label{thm: smallgapoverDelta}
Let $G$ be a graph and $H$ be a $k$-regular graph on $m \geq 1$ vertices. 
Suppose $v$ is a vertex of $G$ for which there are two distinct eigenvalues 
$\lambda, \mu \in \supp_{G}(v)$ such that 
\begin{equation} \label{eq: small-eigenvalue-gap}
	|\lambda - k| - |\mu - k| \ \in \ \{ \sqrt{\Delta}, 2 \sqrt{\Delta}\} 
\end{equation}
for some square-free integer $\Delta$. 
Then $(v,w)$ is not periodic vertex in $G \corona H$, for all $w \in V(H) \cup \{0\}$.
	
\begin{proof}
By Lemma \ref{lem: periodicinbase}, it suffices to show $(v,0)$ is not periodic. 
Suppose towards contradiction that $(v,0)$ is periodic. 
By Lemma \ref{lem: coronaperiodic}, there exists a square-free integer $\Delta$ 
such that for each eigenvalue $\lambda$ in the support of $v$,
both $\lambda - k$ and $\sqrt{(\lambda - k)^2 + 4m}$ are integer multiples of $\sqrt{\Delta}$.
We define
\begin{equation}
	\sigma := \frac{1}{\sqrt{\Delta}} \min\left\{\Bigr| |\lambda_{1} - k| - |\lambda_{2} - k| \Bigr| \ : \
	 \lambda_{1}, \lambda_{2} \in \supp_G(v) \right\}. 
\end{equation}
Let $\lambda$ and $\mu$ be the eigenvalues in the support of $v$ which achieve the above minimum.
Now, define
\begin{equation}
	n_\lambda := \frac{|\lambda-k|}{\sqrt{\Delta}}, \qquad n_\mu := \frac{|\mu-k|}{\sqrt{\Delta}},
\end{equation}
and suppose that $\sigma := n_\lambda - n_\mu$.

By the assumption in \eqref{eq: small-eigenvalue-gap}, we have two cases to consider: $\sigma = 1$ and $\sigma = 2$.

By Lemma \ref{lem: coronaperiodic}, both $n_\lambda^2 + 4m/\Delta$ and
$n_\mu^2 + 4m/\Delta$ are squares.  If we let 
\begin{equation}
	p := \sqrt{n_\mu^2 + \frac{4m}{\Delta}} \quad \text{and} \quad
	q:=\sqrt{n_\lambda^2 + \frac{4m}{\Delta}},
\end{equation}
then
\begin{equation}
q+p > n_\lambda+n_\mu = 2n_\mu+\sigma
\quad \text{and}\quad
q^2-p^2 = (2n_\mu+\sigma)\sigma.
\end{equation}
Hence $q-p < \sigma$ which is impossible when $\sigma=1$.

When $\sigma=2$, $p$ and $q$ have the same parity which contradicts $0< q-p<2$.
\end{proof}
\end{thm}

\begin{cor}  \label{cor: smallgapnotperiodic}
Suppose $G$ is a graph and $H$ is a $k$-regular graph on $m \geq 1$ vertices. 
Let $v \in V(G)$ and suppose there are eigenvalues $\lambda, \mu \in \supp_G(v)$ such that
$0 < |\lambda-k| - |\mu-k| < 3$. Then $(v,w)$ is not periodic in $G \corona H$ for every
$w \in V(H) \cup \{0\}$.

\begin{proof}
Suppose towards contradiction that $(v,0)$ is periodic. Then by Lemma \ref{lem: coronaperiodic},
we have that $\lambda - k$ and $\mu - k$ are integer multiples of $\sqrt{\Delta}$, for some
square-free integer $\Delta$.
Since $0 < |\lambda - k| - |\mu - k| < 3$, we have
\begin{equation}
	|\lambda - k| - |\mu - k| 
	\in \{\sqrt{1}, \sqrt{2}, \sqrt{3}, 2\sqrt{1}, \sqrt{5}, \sqrt{6}, \sqrt{7}, 2\sqrt{2}\}.
\end{equation}
Note that Theorem \ref{thm: smallgapoverDelta} applies to all of these cases.
\end{proof}
\end{cor}

We apply our machinery above to show that the corona products of certain distance-regular graphs
with an arbitrary regular graph have no perfect state transfer. In particular, we show this for some families
of distance-regular graphs which {\em have} perfect state transfer (see Coutinho \etal \cite{cggv15}).
These examples confirm that perfect state transfer is highly sensitive to perturbations.

\begin{cor} \label{cor: drg-corona-no-periodic}
Let $G$ be a graph from one of the following families:
\begin{itemize}
	\item $d$-cubes $Q_d$, for $d \ge 2$.
	\item Cocktail party graphs $\overline{nK_2}$, for $n \ge 2$.
	\item Halved $2d$-cubes $\frac{1}{2}Q_{2d}$, for $d \ge 1$.
	\ignore{
	\item Coset graph of the once shortened and once truncated binary Golay code: \\
	$Spec = \{21, 9, 5, 1, -3, -7, -11\}$.
	\item Coset graph of the shortened binary Golay code:
	$Spec = \{22, 8, 6, 0, -2, -8, -10\}$.
	\item Double coset graph of binary Golay code:
	$Spec = \{23, 9, 7, 1, -1, -7, -9, -23\}$.
	}
\end{itemize}
Then the corona product $G \corona H$, where $H$ is a regular graph, has no perfect state transfer.

\begin{proof}
The spectra of these graphs are well known (see Brouwer \etal \cite{bcn89}) and are given by:
\begin{itemize}
	\item $Spec(Q_d) = \{d - 2\ell : 0 \le \ell \le d\}$.
	\item $Spec(\overline{nK_2}) = \{ 2n - 2, 0, -2 \}$.
	\item $Spec(\frac{1}{2}Q_{2d}) = \left\{ \binom{2d}{2} - 2\ell(2d - \ell) : 0 \le \ell \le d \right\}$.
\end{itemize}
Since these graphs are distance-regular, every eigenvalue is in the support of every vertex.
Moreover, the eigenvalues satisfy the conditions of Corollary \ref{cor: smallgapnotperiodic}.
In particular, $2-d$ and $-d$ are always eigenvalues of the $d$-cube and the halved $2d$-cube.
Therefore, the corona of these graphs with an arbitrary regular graph do not have periodic vertices.
Hence, by Lemma \ref{lem: periodic-necessary-condition}, they do not have perfect state transfer.
\end{proof}
\end{cor}


\subsection{Corona with the Complete Graph}

The main result of this section is that there is no perfect state transfer on $G \corona K_m$ when $G$ is connected 
on at least two vertices. 
This is achieved by bounding the number of vertices in $H$ for the corona product $G \corona H$ to have a periodic vertex.

\begin{lem} \label{lem: boundmk}
Suppose that $G$ is a graph and $H$ is a $k$-regular graph on $m \ge 1$ vertices.
If $(v,0)$ is periodic in $G \corona H$, then for all $\lambda \in \supp_G(v)$ we have
\begin{equation} \label{eqn: size-lower-bound}
m \ \geq \ |\lambda - k| + 1. 
\end{equation}

\begin{proof}
Suppose that $(v,0)$ is periodic in $G \corona H$. 
By Lemma \ref{lem: coronaperiodic}, 
there exists a positive square-free integer $\Delta$ such that for all $\lambda \in \supp_{G}(v)$, 
both $(\lambda - k)^2/\Delta$ and $((\lambda - k)^2 + 4m)/\Delta$ are squares. 
Moreover, $\Delta$ divides $2m$, and therefore, these squares have the same parity.
This yields the bound
\begin{equation}
	\frac{4m}{\Delta} \geq \left( \frac{|\lambda - k|}{\sqrt{\Delta}} + 2 \right)^2- 
	\left(\frac{|\lambda - k|}{\sqrt{\Delta}}\right)^2 
	= 4 \left( \frac{|\lambda - k|}{\sqrt{\Delta}} + 1\right).
\end{equation}
After rearranging, we see that
\begin{equation}
	m \ \ge \ |\lambda - k| \sqrt{\Delta} + \Delta \ \ge \ |\lambda - k| + 1,
\end{equation}
since $\Delta \ge 1$.
\end{proof}
\end{lem}

\begin{thm}
Let $G$ be a connected 
graph on at least 2 vertices. 
Then for all $m \ge 1$, the corona product $G \corona K_m$ has no periodic vertices,
and, therefore, has no perfect state transfer.

\begin{proof}
Suppose towards contradiction that the vertex $(v,w)$ of $G \corona K_m$ is periodic. 
If the eigenvalue support of $v$ in $G$ contains a negative eigenvalue $\lambda$, 
then $\lambda - (m-1) < 0$. Thus, Lemma \ref{lem: boundmk} yields
\begin{equation}
	m \ \geq \ |\lambda - (m-1)| + 1 \ = \ -\lambda + (m-1) + 1 \ > \ m.
\end{equation}
So, it suffices to show that the support of $v$ contains a negative eigenvalue. 

Since $G$ has no loops, we have
\begin{equation}
 \ebra{v}A\eket{v} = \sum_{\lambda \in \Sp{G}} \lambda \ebra{v}\eigenA{\lambda}\eket{v}  =0.
\end{equation}
If $v$ has no negative eigenvalue in its support, then  $\eigenA{\lambda}\eket{v}  =0$ for all $\lambda \neq 0$.
In this case, $A\eket{v}$ is the zero vector and $G$ is not connected.
\end{proof}
\end{thm}


\subsection{Corona with the Empty Graph}

In this section, we show that the corona of any graph with $\emptyg{m}$, 
where $m$ is one or a prime number, has no perfect state transfer.
We will need the following spectral characterization of when a vertex is conical in a star graph.

\begin{lem} \label{lem: starsupport}
Let $v$ be a vertex in a connected graph $G$.
If $\supp_G(v) = \{\pm \lambda\}$, for some $\lambda > 0$, then 
$\lambda^2 \in \mathbb Z$ and $v$ is the conical vertex of the star graph $K_{1,\lambda^2}$.

\begin{proof}
Suppose that $\supp_G(v) = \{\pm \lambda\}$. If $A = \sum_{\theta} \theta \eigenA{\theta}$ is the 
spectral decomposition of the adjacency matrix of $G$, 
then $(\eigenA{\lambda} + \eigenA{-\lambda})\eket{v} = \eket{v}$.
Moreover, $A^{2}\eket{v} = \lambda^{2}\eket{v}$.
Since $A^2$ has integer entries, we observe that $\lambda^2 \in \mathbb Z$. 
This means that every walk of length $2$ starting from $v$ must return to $v$, i.e. that every neighbor 
of $v$ has degree 1. The number of closed walks of length two from $v$ is then exactly the degree of $v$, 
and thus we see that $v$ is the center vertex of $K_{1,\lambda^2}$.
\end{proof}
\end{lem}

\begin{lem} \label{lem: noPSTclaw}
$K_{1,n} \corona \emptyg{m}$ has no perfect state transfer, for every $n,m \ge 1$.

\begin{proof}
The case for $n = 1$ is a result of Fan and Godsil \cite{fg13}.
For $n = 2$, that is, $P_{3} \corona \emptyg{m}$, we may apply Theorems \ref{thm: godsil-automorphism} 
and \ref{thm: coutinho-pst-conditions}
and Lemma \ref{lem: coronaperiodic} to rule out perfect state transfer.
For $n > 2$ and $m \neq 2$, we may apply Theorems \ref{thm: godsil-automorphism} 
and \ref{thm: coutinho-pst-conditions} to show no perfect state transfer exists.

So, we consider $K_{1,n} \corona \emptyg{2}$ where $n > 2$, and assume that it has perfect state transfer 
between vertices $(v,1)$ and $(v,2)$ where $v$ is a vertex of $K_{1,n}$. 
Note that perfect state transfer between other pairs of vertices are ruled out by 
Theorems \ref{thm: godsil-automorphism} and \ref{thm: coutinho-pst-conditions}.
Thus, $(v,0)$ is a periodic vertex.
By Lemma \ref{lem: coronaperiodic}, there is a square-free integer $\Delta$ so that 
for each eigenvalue $\lambda \in \supp(v)$, both $\lambda$ and $\sqrt{\lambda^{2} + 8}$ 
are integer multiples of $\sqrt{\Delta}$. Moreover, $\Delta$ divides $4$ and thus
$\Delta \in \{1,2\}$.

If $\Delta = 1$, both $\lambda$ and $\sqrt{\lambda^{2} + 8}$ are integers 
for each eigenvalue $\lambda \in \supp(v)$. 
Thus, $\supp(v) = \{\pm 1\}$. 
By Lemma \ref{lem: starsupport}, $v$ is a conical vertex of $K_{1,1}$.
But, $K_{1,1} \corona \emptyg{2}$ has no perfect state transfer (see \cite{fg13}).

If $\Delta = 2$, both $\lambda$ and $\sqrt{\lambda^{2} + 8}$ are integer multiples of $\sqrt{2}$
for each eigenvalue $\lambda \in \supp(v)$. 
Suppose $\lambda = \ell\sqrt{2}$ for some integer $\ell$. 
Then, $\ell^{2} + 4$ is a square which implies $\ell = 0$. 
Thus, $\supp(v) = \{0\}$. But, this implies $v$ is an isolated vertex, which is a contradiction.
\end{proof}
\end{lem}

\begin{thm}\label{thm: noPSThairs}
If $G$ is a connected graph on at least two vertices and $m$ is either $1$ or a prime number,
then $G \corona \emptyg{m}$ has no perfect state transfer.

\begin{proof}
By Lemma \ref{lem: noPSTclaw}, we may assume $G$ is connected and is not a star.
By Lemma \ref{lem: periodic-necessary-condition}, it suffices to show 
$G \corona \emptyg{m}$ has no periodic vertices for $m \ge 1$.
Suppose for contradiction that vertex $(v,w)$ is periodic in $G \corona \overline{K}_{m}$.
Then by Lemma \ref{lem: periodicinbase}, $(v,0)$ is periodic.
Since $G \ne K_{1,n}$ is connected, by Lemma \ref{lem: starsupport}, the eigenvalue support of $v$ in $G$ 
contains eigenvalues $\lambda$ and $\mu$ such that $|\lambda| \ne |\mu|$. 
Assume $|\mu| < |\lambda|$. 
By Lemma \ref{lem: coronaperiodic}, there must exist a square-free integer 
$\Delta$ dividing $2m$ such that $\lambda$, $\mu$, $\sqrt{\lambda^2 + 4m}$, and $\sqrt{\mu^2 + 4m}$ 
are all integer multiples of $\sqrt{\Delta}$. 

We define the integers 
\begin{equation}
n_\lambda := \frac{|\lambda|}{\sqrt{\Delta}}, 
\ \hspace{0.2in} \
n_\mu := \frac{|\mu|}{\sqrt{\Delta}},
\ \hspace{0.2in} \
\ell := \frac{2m}{\Delta}.
\end{equation} 
Since $\sqrt{\mu^2 + 4m}$ and $\sqrt{\lambda^2 + 4m}$ are integer multiples of $\sqrt{\Delta}$,
both $n_\mu^2 + 2\ell$ and $n_\lambda^2 + 2\ell$ are perfect squares.
So, let $N_{\mu}^{2} := n_\mu^2 + 2\ell$ and $N_{\lambda}^{2} := n_\lambda^2 + 2\ell$. 
Therefore, we have four perfect squares: 
$n_{\mu}^{2}$, $N_{\mu}^{2}$, $n_{\lambda}^{2}$, and $N_{\lambda}^{2}$.
Since $|\mu| < |\lambda|$, we know $n_{\mu} < n_{\lambda}$. 

If $m=1$ or $m=2$, then $\ell$ equals one, two or four and $2\ell$ equals two, four or eight.
First note that a difference of two squares is never two or four.
For the last case, the only way for $N_{\mu}^{2}-n_{\mu}^2 = N_{\lambda}^2-n_{\lambda}^2=8$ is 
$N_{\lambda}=N_{\mu}=3$ and $n_{\mu}=n_{\lambda}=1$, contradicting $n_{\lambda}>n_{\mu}$.

Now, consider when $m$ is an odd prime. Since a squarefree integer $\Delta$ divides $2m$, 
we have $\ell$ equals one, two, $m$ or $2m$. Or, equivalently, $2\ell$ equals two, four, $2m$ or $4m$.
We may rule out the first two cases as before.
Since $N_{\mu}^{2} - n_{\mu}^{2}$ is even, $N_{\mu}$ and $n_{\mu}$ share the same parity. 
If $2\ell = 2m$, then our odd prime $m$ can be written as 
$m = (N_{\mu} - n_{\mu}) (N_{\mu}+n_{\mu})/2$.
This shows $N_{\mu}-n_{\mu}=1$ or that $N_{\mu}$ and $n_{\mu}$ differ in parity, 
a contradiction.
If $2\ell = 4m$, then 
$m = \frac{1}{2}(N_{\mu} - n_{\mu}) \frac{1}{2}(N_{\mu}+n_{\mu})$.
This shows $N_{\mu} - n_{\mu} = 2$, and furthermore, $m = n_{\mu} + 1$.
But using the other difference of squares, we also have $m = n_{\lambda} + 1$.
This contradicts $n_{\lambda} > n_{\mu}$.
\end{proof}
\end{thm}


\subsection{Multigraphs}

Fan and Godsil \cite{fg13} proved that the double star $K_{2} \corona \emptyg{m}$ has no perfect state transfer 
for all $m$.
Here, we show that $C_{2} \corona \emptyg{m}$ {\em has} perfect state transfer for some $m$ where
$C_{2}$ is the {\em digon} (a multigraph on two vertices which are connected by two parallel edges).

This result follows from the next proposition.
For a real number $\alpha$, let $K_{2}(\alpha)$ denote the clique on two vertices whose edge has weight $\alpha$.
In what follows, we consider the weighted path $K_{2}(\alpha) \corona K_{1}$.

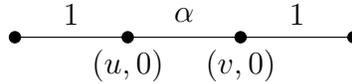
\begin{figure}[H]
	\centering
	\topline{.15in}
	\begin{tikzpicture}[scale=0.75]
		\filldraw (0,0) circle[radius=0.1] -- (2,0) circle[radius=0.1] -- (4,0) circle[radius=0.1] -- (6,0) circle[radius=0.1];
		\draw (1,0.4) node{$1$} ++(2,0) node{$\alpha$} ++(2,0) node{$1$};
		\draw (2,-0.5) node{$(u,0)$} ++(2,0) node{$(v,0)$};
	\end{tikzpicture}
	\caption{The weighted path $K_{2}(\alpha) \corona K_{1}$.} \label{fig: 1alpha1}
	\botline{.15in}
\end{figure}

\begin{prop} \label{prop: weighted-p4}

In $K_{2}(\alpha) \corona K_{1}$, perfect state transfer occurs between the vertices of $K_{2}$
if and only if 
\begin{equation}
\alpha = \frac{2(2s+1)}{\sqrt{(2\ell)^2 - (2s+1)^2}}
\qquad \text{ for some integers $\ell>s\geq 0$.}
\end{equation}
In this case, perfect state transfer occurs at time
$t = \frac{\pi}{2} \sqrt{(2\ell)^2 - (2s+1)^2}$.

\begin{proof}
Let $u, v$ denote the vertices of $K_2$. 
The eigenvalues of $K_2(\alpha)$ are $\pm \alpha$ with corresponding eigenprojectors 
$\eigenA{\pm\alpha}$ which satisfy
$\ebra{u} \eigenA{\pm\alpha} \eket{v} = \pm \frac{1}{2}$.
By Proposition \ref{prop: adj-transition-corona}, the transition element between $(u,0)$ and $(v,0)$ is given by
\begin{equation}
	\ebra{(u,0)} e^{-\ii tA(K_2(\alpha)\corona K_1)} \eket{(v,0)}
		= -i \sin\left(\frac{t}{2}\alpha\right)\cos\left(\frac{t}{2}\Damp_\alpha\right)
			-i\frac{\alpha}{\Damp_\alpha} \cos\left(\frac{t}{2}\alpha\right)\sin\left(\frac{t}{2}\Damp_\alpha\right),
\end{equation}
where $\Damp_\alpha := \sqrt{\alpha^2 + 4}$. 
Since $|\alpha/\Damp_\alpha|<1$, there is perfect state transfer at time $t$ if and only if 
\begin{equation}
	\left| \sin\left(\frac{t}{2}\alpha\right) \cos\left( \frac{ t}{2}\Damp_\alpha\right) \right| = 1.
\end{equation}
Equivalently, this gives the conditions
\begin{align}
\frac{t}{2} \sqrt{\alpha^2 + 4} = \ell\pi,
\qquad \text{and} \qquad
\frac{t}{2} \alpha = \left( s + \frac{1}{2} \right) \pi,
\end{align}
for some integers $\ell>s\geq0$.
Hence perfect state transfer occurs at time $t$ if and only if
\begin{equation}
\alpha = \frac{2(2s+1)}{\sqrt{(2\ell)^2 - (2s+1)^2}}.
\end{equation}
\end{proof}
\end{prop}

\begin{figure}[h]
\centering
\topline{.15in}
\begin{tikzpicture}[scale = 0.85]
	\foreach \x in {0,1,2,9,10} {
		\filldraw (2,0) -- ++(-2, 2.5 - 0.5*\x) circle[radius=0.1];
		\filldraw (4,0) -- ++(2, 2.5 - 0.5*\x) circle[radius=0.1];
	}
	\draw[dashed] (0,0.5) -- (0,-1);
	\draw[dashed] (6,0.5) -- (6,-1);
	\draw (2,0) to[bend right] ++(2,0) to[bend right] ++(-2,0);
	\filldraw[fill=white] (2,0) circle[radius=0.1] ++(0,-0.4) node{$u$};
	\filldraw[fill=gray] (4,0) circle[radius=0.1] ++(0,-0.4) node{$v$};
\end{tikzpicture} \quad \quad \quad \quad
\begin{tikzpicture}[scale = 0.85]
   	\draw (0,0) -- ++(1.5,1.5) -- ++(1.5,-1.5) -- ++(-1.5,-1.5) -- ++(-1.5,1.5);
	\foreach \x in {0,1,2,9,10} {
		\filldraw (0,0) -- ++(-1, 1.5 - 0.3*\x) circle[radius=0.1]; 
		\filldraw (1.5,1.5) -- ++(-1.5 + 0.3*\x, 1) circle[radius = 0.1]; 
		\filldraw (3,0) -- ++(1, 1.5 - 0.3*\x) circle[radius=0.1]; 
		\filldraw (1.5,-1.5) -- ++(-1.5 + 0.3*\x, -1) circle[radius = 0.1]; 
	} 
	\draw[dashed] (-1, 0.5) -- ++(0, -1); 
	\draw[dashed] (1, 2.5) -- ++(1, 0); 
	\draw[dashed] (4, 0.5) -- ++(0, -1); 
	\draw[dashed] (1, -2.5) -- ++(1, 0); 
	\filldraw[fill=white] (0,0) circle[radius=0.1] ++(0,-0.65) node{$u_1$}; 
	\filldraw[fill=white] (3,0) circle[radius=0.1] ++(0,-0.65) node{$u_2$}; 
	\filldraw[fill=gray] (1.5,1.5) circle[radius=0.1] ++(0.65,0) node{$v_1$}; 
	\filldraw[fill=gray] (1.5,-1.5) circle[radius=0.1] ++(0.65,0) node{$v_2$}; 
\end{tikzpicture} 

\caption{
For a positive integer $r$,
$C_2 \corona \emptyg{4r^2-1}$ has perfect state transfer between $u$ and $v$,
and
$C_{4} \corona \emptyg{4r^{2} - 1}$ has {\em real} perfect state transfer between
$\{u_1,u_2\}$ and $\{v_1, v_2\}$.}
\label{fig:coronaempty}
\botline{.15in}
\end{figure}
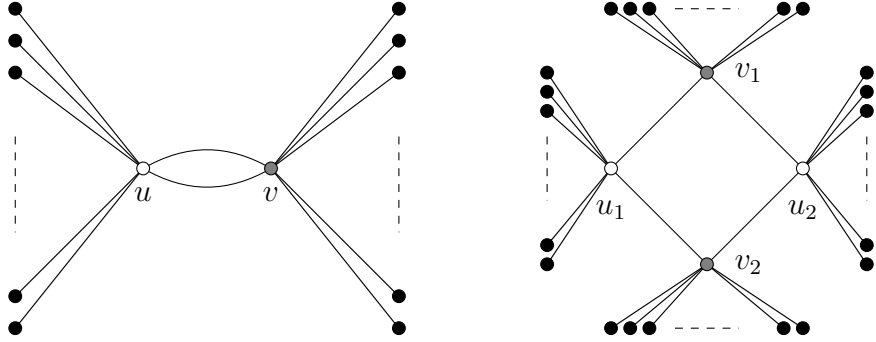

In a graph $G = (V,E)$, for a subset $U \subseteq V$ of vertices, we denote
the real uniform superposition of vertices in $U$
as $\eket{U} = |U|^{-1/2}\sum_{v \in U} \eket{v}$. 
We say $G$ has {\em real} perfect state transfer between two subsets 
$U_{1},U_{2} \subseteq V$ if $|\ebra{U_{1}}e^{-\ii tA(G)}\eket{U_{2}}| = 1$.

\begin{thm} \label{thm: digon-double-star}
For each positive integer $r$, the following hold.
\begin{enumerate}[i)]
\item 
$C_{2} \corona \emptyg{4r^{2}-1}$ has perfect state transfer
between the vertices of $C_{2}$ at time $\pi/2$. 

\item 
$C_{4} \corona \emptyg{4r^{2}-1}$ has real perfect state transfer
between its two antipodal pairs of vertices
at time $\pi/2$.
\end{enumerate}

\begin{proof} 
We apply Proposition \ref{prop: weighted-p4} with $s=0$ and $\ell=r$ for a positive integer $r$.
This shows that $P_{4}(r) := K_{2}(2/\sqrt{4r^{2}-1}) \corona K_1$
has perfect state transfer between the inner two vertices at time $\tau = \frac{\pi}{2}\sqrt{4r^{2}-1}$.
Hence there exists a unimodular complex number $\gamma$ such that
$e^{-i\tau A_(P_4(r))} e_{(u,0)}= \gamma e_{(v,0)}$.
If follows from 
\begin{equation}
A(P_4(r)) \left( e^{-i\tau A(P_4(r))} e_{(u,0)} \right) = e^{-i\tau A(P_4(r))}A(P_4(r)) e_{(u,0)} = \gamma A(P_4(r)) e_{(v,0)}
\end{equation}
 that perfect state transfer also occurs between the antipodal vertices.

By multiplying all weights of this graph by $\sqrt{4r^{2}-1}$, 
we obtain a weighted graph $\mathcal{G}_{r}$ whose adjacency matrix is
$A(\mathcal{G}_{r}) = \sqrt{4r^{2}-1} A(P_{4}(r))$;
see Figure \ref{fig: weighted-p4}.
Note $\mathcal{G}_{r}$ has perfect state transfer at time $\pi/2$ between the two inner vertices.

\begin{figure}[H] \label{fig: weightedpathpst}
	\centering
	\topline{.15in}
	\begin{tikzpicture}
		\filldraw (0,0) circle[radius=0.1] -- (2,0) circle[radius=0.1] -- (4,0) circle[radius=0.1] -- (6,0) circle[radius=0.1];
		\draw (1,0.3) node{$\sqrt{4r^2 - 1}$} ++(2,0) node{$2$} ++(2,0) node{$\sqrt{4r^2 - 1}$};
		\draw (2,-0.5) node{$u$} ++(2,0) node{$v$};
		\draw (0,-0.5) node{$a$} ++(6,0) node{$b$};
	\end{tikzpicture}
	\caption{A family of weighted paths $\mathcal{G}_{r}$, for positive integer $r$, with perfect state transfer 
		(between $u$ and $v$, and between $a$ and $b$). 
	}
	\label{fig: weighted-p4}
	\botline{.15in}
\end{figure}
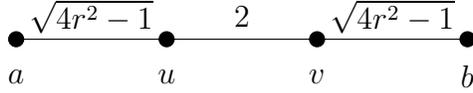

\par\noindent
Moreover, $\mathcal{G}_{r}$ is the quotient of both $C_{2} \corona \emptyg{4r^{2} - 1}$ 
and $C_{4} \corona \emptyg{4r^{2} - 1}$ under natural equitable partitions.
Since perfect state transfer is closed under taking quotient and lifting
(see Bachman \etal \cite{bfflott12}, Theorem 1), this completes both claims.
\end{proof}
\end{thm}



Consider the family of weighted paths $\mathcal{G}_{r}$ given in Figure \ref{fig: weighted-p4}.
The weighted path $\mathcal{G}_{1}$ corresponds to the distance quotient of the cube $Q_{3}$ 
(see also Figure \ref{fig: q3-quotient-lift}). It is known that the distance quotients of the cube $Q_{n}$ 
provide a family of weighted paths with perfect state transfer. The weighting schemes of these weighted paths
are unimodal (single peaked), where the edge weight between the $k$th and $(k+1)$th vertices is
$\sqrt{(k+1)(n-k)}$ for $k=0,\ldots,n-1$.
In contrast, the weighted paths $\mathcal{G}_{r}$, for $r > 1$, exhibit weighting schemes that are not unimodal.

\newcommand*{\nodeRadius}{0.08}
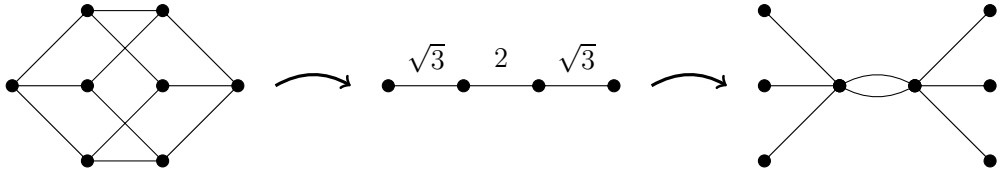
\begin{figure}[H]
	\centering
	\topline{.15in}
	\begin{tikzpicture}[scale=1.0]
		\filldraw (0,0) circle[radius=\nodeRadius];
		\foreach \x in {-1,0,1} {
			\filldraw (0,0) -- (1, \x) circle[radius=\nodeRadius];
		}
		\filldraw (3,0) circle[radius=\nodeRadius];
		\foreach \x in {-1,0,1} {
			\filldraw (3,0) -- (2, \x) circle[radius=\nodeRadius];
		}
		\draw (1,1) -- (2,1) -- (1,0) -- (2,-1) -- (1,-1) -- (2,0) -- (1,1);
		
		\draw[very thick, ->] (3.5,0) to[bend left] ++(1,0);
		
		\filldraw (5,0) circle[radius=\nodeRadius] -- ++(1,0) circle[radius=\nodeRadius]
		 -- ++(1,0) circle[radius=\nodeRadius] -- ++(1,0) circle[radius=\nodeRadius];
		 \draw (5.5, 0.375) node[scale=0.9]{$\sqrt{3}$} ++ (1,0) node[scale=0.9]{$2$} ++(1,0) node[scale=0.9]{$\sqrt{3}$};
		 
		 \draw[very thick, ->] (8.5,0) to[bend left] ++(1,0);
		 
		 \draw (11, 0) to[bend right] ++(1,0) to[bend right] ++(-1,0);
		 \foreach \x in {-1,0,1} {
		 	\filldraw (10, \x) circle[radius=\nodeRadius] -- (11,0) circle[radius=\nodeRadius];
		 }
		 \foreach \x in {-1,0,1} {
		 	\filldraw (13, \x) circle[radius=\nodeRadius] -- (12,0) circle[radius=\nodeRadius];
		 }
	\end{tikzpicture}
	\caption{Example: a quotient and a lift of the hypercube $Q_3$.}
	\label{fig: q3-quotient-lift}
	\botline{.15in}
\end{figure}


\section{Pretty Good State Transfer}

\subsection{Barbell Graphs} \label{sub: k2coronaclique}

We consider the family of barbell graphs (see Ghosh \etal \cite{gbs08}) and show
that they exhibit pretty good state transfer.

\begin{figure}[H]
\centering
\newcommand*{\cliqueRadius}{1}
\newcommand*{\cliqueVertices}{7}
 \begin{tikzpicture}[scale = 1.25]
	\foreach \s in {0,180} 
	{
		\pgfmathsetmacro{\cv}{\cliqueVertices + 1}
		\foreach \t in {1,...,\cv}
		{
			\pgfmathsetmacro{\angle}{360*\t/(\cv) + 180}
			\filldraw {(\s:2) + (\angle + \s: \cliqueRadius)} circle[radius=0.075];
			\foreach \x in {1,...,\cv}
			{
				\pgfmathsetmacro{\anglex}{360*\x/(\cv) +180}
				\draw ($(\s:2) + (\angle + \s: \cliqueRadius)$) -- ($(\s:2) + (\anglex + \s: \cliqueRadius)$);
			}
			\draw (0,0) -- (\s:2-\cliqueRadius);
			\filldraw[color=black, fill=white] (\s:2-\cliqueRadius) circle[radius=0.075];
		}
	}
\end{tikzpicture}
\caption{$K_2 \corona K_{\cliqueVertices}$ has pretty good state transfer between vertices marked white.}
\botline{.15in}
\end{figure}
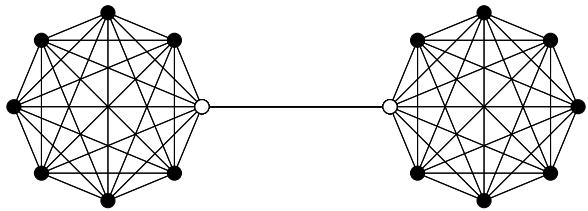

\begin{thm}\label{thm: clique pretty good state transfer}
$K_2 \corona K_m$ has pretty good state transfer between the vertices of $K_2$ for all $m \geq 1$.

\begin{proof}
Let $u$ and $v$ denote the vertices of $K_2$. 
The eigenvalues of $K_2$ are $\pm 1$, with $\ebra{u}\eigenA{\pm}{K_2}\eket{v} = \pm 1/2$. 
Let $m \ge 1$ be given. We show there is pretty good state transfer between 
$(u,0)$ and $(v,0)$ in $K_{2} \corona K_{m}$.

In Proposition \ref{prop: adj-transition-corona}, we have $\Damp_{\pm} = \sqrt{((m-1) \mp 1)^{2} + 4m}$
which are irrational for $m \ge 1$.
If $\Damp_{\pm}$ are linearly independent over $\Q$,
by Kronecker's Theorem, there are integers $\ell$ and $r_{\pm}$ such that
\begin{equation} \label{eq: barbells-independent-approximation}
	\ell \Damp_{\pm}/2  - r_{\pm} \ \approx \ - \Damp_{\pm}/4.
\end{equation}
This shows that
\begin{equation}
	(2\ell + 1)\pi \Damp_{\pm}/2 \ \approx \ 2\pi r_{\pm}.
\end{equation}
So, if we let $t = (2\ell + 1)\pi$, then
\begin{equation}
|e^{-\ii((m-1) \pm 1)t/2} \cos\left(t\Damp_{\pm}/2\right)| \ \approx \ 1.
\end{equation}
By Proposition \ref{prop: adj-transition-corona}, we get
$|\ebra{(u,0)} e^{-\ii tA(G \corona H)} \eket{(v,0)}| \approx 1$.

On the other hand, if $\Damp_{\pm}$ are linearly dependent over $\Q$, then there are integers
$a_{\pm}$ and a square-free integer $\Theta$ so that $\Damp_{\pm} = a_{\pm}\sqrt{\Theta}$.
Since $m^{2}+4$ is not a perfect square, $\sqrt{\Theta}$ is irrational.
By Kronecker's Theorem, there are integers $\ell$ and $q$ such that
\begin{equation} \label{eq: barbells-dependent-approximation}
	\ell \sqrt{\Theta}/2 - q \ \approx \ -\sqrt{\Theta}/4.
\end{equation}
This shows that
\begin{equation}
	(2\ell + 1)\pi \Damp_{\pm}/2 \ \approx \ 2\pi a_{\pm}q.
\end{equation}
Similarly, if we let $t = (2\ell + 1)\pi$, by Proposition \ref{prop: adj-transition-corona}, we get
$|\ebra{(u,0)} e^{-\ii tA(G \corona H)} \eket{(v,0)}| \approx 1$.
\end{proof}
\end{thm}


\subsection{Hairy Cliques}

By Theorem \ref{thm: noPSThairs}, we know that there is no perfect state transfer between 
$(v,1)$ and $(v,2)$ in $G \corona \emptyg{2}$ unless $v$ is isolated in $G$. 
We show a stronger result that there is no pretty good state transfer when $G$ is the complete graph.

Let $K_n + I$ denote the graph whose adjacency matrix is the all-one matrix $J_{n}$. 
That is, $K_n + I$ is obtained from the complete graph by adding self-loops to each vertex. 
Although the quantum walks on $K_n$ and on $K_n + I$ are equivalent up to phase factors,
we show that in contrast to Theorem \ref{thm: noPSThairs} with $K_{n} \corona \emptyg{2}$, 
the corona $(K_{n} + I) \corona \emptyg{2}$ has pretty good state transfer between 
the vertices of $\emptyg{2}$,

\newcommand*{\kRadius}{2}
\newcommand*{\hairRadius}{3}
\newcommand*{\vertexNumber}{8}
\newcommand*{\hairAngle}{10}
\newcommand*{\cliquescale}{0.6}
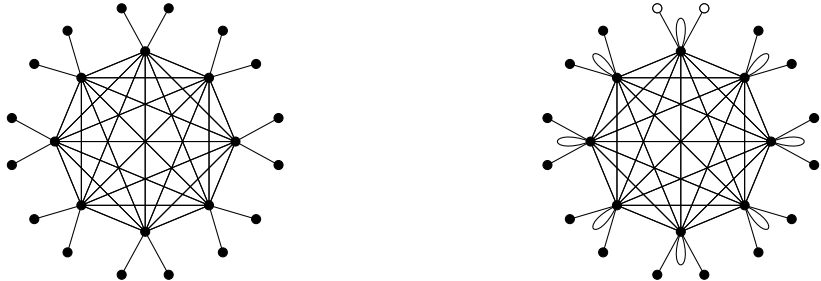
\begin{figure}[H]
\centering
\begin{tikzpicture}[scale = \cliquescale]
	\foreach \t in {1, ..., \vertexNumber }
	{
		\pgfmathsetmacro{\anglet}{(360 / \vertexNumber)*\t + 90}
		\filldraw (\anglet : \kRadius) circle[radius=0.1];
		\foreach \s in {-\hairAngle,\hairAngle}
		{
			\filldraw (\anglet + \s: \hairRadius) circle[radius=0.1];
			\draw (\anglet : \kRadius) -- (\anglet + \s: \hairRadius); 
		}
		\foreach \s in {1, ..., \vertexNumber}
		{
			\pgfmathsetmacro{\angles}{(360/\vertexNumber)* \s + 90}
			\draw (\anglet: \kRadius) -- (\angles : \kRadius);
		}
	}
\end{tikzpicture}
\qquad \qquad \qquad \qquad
\begin{tikzpicture}[scale = \cliquescale]
	\foreach \t in {1, ..., \vertexNumber }
	{
		\pgfmathsetmacro{\anglet}{(360 / \vertexNumber)*\t + 90}
		\filldraw (\anglet : \kRadius) circle[radius=0.1];
		\foreach \s in {-\hairAngle,\hairAngle}
		{
			\filldraw (\anglet + \s: \hairRadius) circle[radius=0.1];
			\draw (\anglet : \kRadius) -- (\anglet + \s: \hairRadius); 
		}
		\foreach \s in {1, ..., \vertexNumber}
		{
			\pgfmathsetmacro{\angles}{(360/\vertexNumber)* \s + 90}
			\draw (\anglet: \kRadius) -- (\angles : \kRadius);
		}
		\draw (\anglet: \kRadius) .. 
			controls (\anglet + 2*\hairAngle/3: \hairRadius) and (\anglet -2*\hairAngle/3: \hairRadius) .. 
			(\anglet: \kRadius);
	}
	\foreach \s in {-\hairAngle + 90, \hairAngle + 90}
	{
		\filldraw[color=black, fill=white] (\s: \hairRadius ) circle[radius = 0.1];
	}
\end{tikzpicture}
\caption{$K_{\vertexNumber}\corona \emptyg{2}$ has no pretty good state transfer, 
while $(K_{\vertexNumber} + I) \corona \emptyg{2}$ has pretty good state transfer between 
every two vertices of degree $1$ at distance $2$.}
\botline{.15in}
\end{figure}

\begin{thm} \label{thm: selflooppgst1}
For $n \ge 3$, there is no pretty good state transfer in $K_n \corona \emptyg{2}$.

\begin{proof}
By Theorem \ref{thm:pgst-necessary-conditions},
we may rule out pretty good state transfer
from $(v,0)$ to $(v,a)$ for $a=1,2$, and from $(v,a)$ to $(w,b)$ for $a,b=0,1,2$ and $v\neq w$.
\ignore{
from $(v,0)$ to $(w,b)$, for $b=1,2$,
from $(v,1)$ to $(w,b)$, for $v \neq w$ and $b=1,2$,
and
from $(v,0)$ to $(w,0)$, for $v \neq w$.
}

Next, we rule out pretty good state transfer between $(v,1)$ and $(v,2)$ for each vertex $v$ of $K_{n}$.
By Proposition \ref{prop: adj-transition-corona}, we may write
\begin{equation}\label{eq: transition between hairs}
\ebra{(v,1)} e^{-\ii tA(K_{n} \corona \emptyg{2})} \eket{(v,2)} 
	= 
		(\ebra{1}\ket{z}_{0}) (\tbra{z}_{0}\eket{2})
		+ \sum_{\lambda \in \Sp(K_{n})} 
			\sum_{\pm} \frac{e^{-\ii t\lambda_{\pm}}}{\lambda_{\pm}^{2}+2} 
			\ebra{v}\eigenA{\lambda}{K_{n}}\eket{v}
\end{equation}
where $\lambda_{\pm} = \frac{1}{2}(\lambda \pm \sqrt{\lambda^{2}+8})$
and $\ket{z}_{0} = \frac{1}{\sqrt{2}}(\eket{1}-\eket{2})$ is an eigenvector of $0$ in $\emptyg{2}$
orthogonal to the all-one eigenvector.
After some straightforward algebraic manipulation, we get
\begin{equation}
\sum_{\pm} \frac{e^{-\ii t\lambda_{\pm}}}{\lambda_{\pm}^{2}+2} 
	\ = \
	\frac{e^{-\ii\lambda t/2}}{2}
	\left[ \cos\left( \frac{t}{2}\sqrt{\lambda^2 + 8}\right) + 
	\frac{\ii \lambda}{\sqrt{\lambda^2 + 8}} \sin\left( \frac{t}{2} \sqrt{\lambda^2 + 8} \right) \right].
\end{equation}
Therefore, the transition matrix element $\ebra{(v,1)} e^{-\ii tA(K_{n} \corona \emptyg{2})} \eket{(v,2)}$
is equal to
\begin{multline} \label{eq: transition between hairs second part}
-\frac{1}{2} + \frac{1}{2} \sum_{\lambda \in \Sp(K_{n})} e^{-\ii\lambda t / 2} 
	\left[ \cos\left( \frac{t}{2}\sqrt{\lambda^2 + 8}\right) + 
	\frac{\ii \lambda}{\sqrt{\lambda^2 + 8}} \sin\left(\frac{t}{2}\sqrt{\lambda^2 + 8}\right) \right] 
	\ebra{v} \eigenA{\lambda}{K_{n}}\eket{v}.
\end{multline}
Since the spectra of $K_{n}$ is $\{(n-1)^{(1)},(-1)^{(n-1)}\}$, we see that
\begin{eqnarray} \label{eq:hairy-clique}
\lefteqn{\ebra{(v,1)} e^{-\ii tA(K_{n} \corona \emptyg{2})} \eket{(v,2)}} \\
	& = &
	-\frac{1}{2} 
	+ 
	\frac{(n-1)e^{\ii t/2}}{2n} 
	\left[ \cos\frac{3t}{2} - \frac{\ii}{3} \sin\frac{3t}{2} \right] 
	+ 
	\frac{e^{-\ii(n-1)t/2}}{2n} 
	\left[ \cos\frac{\Lambda t}{2} + \frac{\ii(n-1)}{\Lambda} \sin\frac{\Lambda t}{2} \right]
\end{eqnarray}
where $\Lambda = \sqrt{(n-1)^{2} + 8}$. 
By Lemma \ref{lem:pgst-necessity}, a necessary condition for pretty good state transfer between
$(v,1)$ and $(v,2)$ is that 
\begin{equation}
\cos\frac{t}{2} \cos\frac{3t}{2} \approx -1
\end{equation}
which is impossible.
\end{proof}
\end{thm}

\begin{thm} \label{thm: selflooppgst2}
Let $K_{n} + I$ denote the graph obtained by adding self-loop to each vertex of the complete graph $K_{n}$,
where $n \ge 1$. Then for $n \ge 2$, 
there is pretty good state transfer in $(K_n + I) \corona \emptyg{2}$ between vertices $(v,1)$ and $(v,2)$, 
for each vertex $v$ of $K_{n}$, where the vertices of $\emptyg{2}$ are denoted by $1$ and $2$.

\begin{proof}
The spectrum of $K_{n} + I$ is given by $\{n,0\}$ with the same eigenspaces as $K_{n}$. 
Thus, the transition matrix element $\ebra{(v,1)}e^{-\ii tA((K_{n} + I) \corona \emptyg{2})}\eket{(v,2)}$
may be derived from \eqref{eq: transition between hairs second part}
since $\eigenA{\lambda + 1}(K_{n} + I) = \eigenA{\lambda}(K_{n})$.
Also, each vertex of $K_{n} + I$ has the full eigenvalue support.
Similar to \eqref{eq: transition between hairs second part}, 
sufficient conditions for pretty good state transfer are:
\begin{eqnarray} 
	\label{eq: pgst-clique-selfloop1}
	\cos\left( \frac{t}{2}\sqrt{8} \right) & \approx & -1, \\ 
	\label{eq: pgst-clique-selfloop2}
	\cos\left( \frac{t}{2} n \right)\cos\left( \frac{t}{2}\sqrt{n^2 + 8} \right) & \approx & -1.
\end{eqnarray}
Since $\sqrt{n^{2} + 8}$ is irrational for $n \ge 2$, both
$n/\sqrt{2}$ and $\sqrt{(n^{2} + 8)/2}$ are linearly independent over $\Q$.
We consider two cases based on whether $(n^{2} + 8)/2$ is a perfect square or not.

If $(n^{2} + 8)/2$ is not a perfect square, by Kronecker's Theorem,
there are integers $r$, $s$, and $\ell$ so that
\begin{eqnarray}
	\label{eq: kronecker-pgst-clique-selfloop1}
	r \left( \frac{n}{2\sqrt{2}} \right) - s & \approx & - \frac{n}{4\sqrt{2}}, \\
	\label{eq: kronecker-pgst-clique-selfloop2}
	r \left( \frac{1}{2} \sqrt{\frac{n^2+8}{2}} \right) - \ell & \approx & \frac{1}{2} - \frac{1}{4} \sqrt{\frac{n^2+8}{2}}.
\end{eqnarray}
So, if we let $t = (2r + 1)\pi/\sqrt{2}$, then the last equations
are equivalent to
$tn/2 \approx 2\pi s$
and
$t\sqrt{n^{2} + 8}/2 \approx (2\ell + 1)\pi$.
Therefore, \eqref{eq: pgst-clique-selfloop1} and \eqref{eq: pgst-clique-selfloop2} are satisfied
and we have pretty good state transfer.

Otherwise, if $(n^{2} + 8)/2$ is a perfect square, then both $n$ and $(n^{2} + 8)/2$ must be even.
So, let $(n^{2} + 8)/2 = (2x)^{2}$ and $n = 2y$, for some integers $x$ and $y$.
Substituting for $n$, we get $y^{2} + 2 = 2x^{2}$. Reasoning modulo $4$, we see $x$ must be odd,
say, $x = 2q + 1$ for some integer $q$.
First, we apply Kronecker's Theorem to obtain the integers $r,s$ 
to satisfy \eqref{eq: kronecker-pgst-clique-selfloop1},
and then we let $\ell = r(2q + 1) + q$ to obtain
an exact solution to \eqref{eq: kronecker-pgst-clique-selfloop2}.
As before, this satisfies \eqref{eq: pgst-clique-selfloop1} and \eqref{eq: pgst-clique-selfloop2}.
Hence, we have pretty good state transfer.
\end{proof}
\end{thm}

We observe that $(K_1 + I) \corona \emptyg{2}$ is periodic but has no pretty good state transfer.
The eigenvalues of the corona $(K_{1} + I) \corona \emptyg{2}$ are $2$, $0$, and $-1$, which shows periodicity.
However, there is no perfect state transfer, since that would require
$\cos(t/2) \cos(3t/2) = -1$, which is impossible. 
Because there is no perfect state transfer and the graph is periodic, there can be no pretty good state transfer
(see Fan and Godsil \cite{fg13}).

\renewcommand*{\kRadius}{2}
\renewcommand*{\hairRadius}{3}
\renewcommand*{\vertexNumber}{8}
\renewcommand*{\hairAngle}{10}
\renewcommand*{\cliquescale}{0.65}
\newcommand*{\dx}{2.05cm}
\newcommand*{\rotAngleX}{7.25}
\newcommand*{\rotAngleY}{17.5}
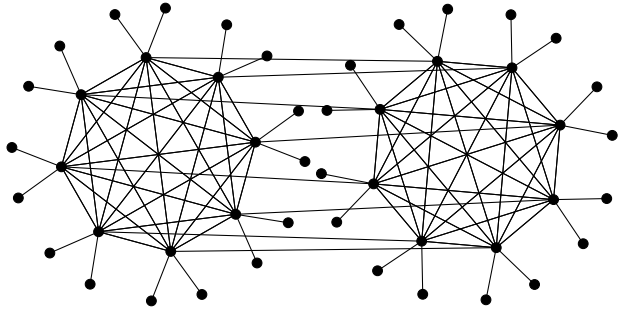
\begin{figure}[H]
\centering
\topline{.15in}
\begin{tikzpicture}[scale = \cliquescale]
	\foreach \t in {1, ..., \vertexNumber }
	{
		\pgfmathsetmacro{\anglet}{(360 / \vertexNumber)*\t + \rotAngleX}
		\coordinate[shift={(-\dx, 0cm)}] (lPoint) at (\anglet : \kRadius);
		\filldraw (lPoint) circle[radius=0.1];
		\foreach \s in {-\hairAngle,\hairAngle}
		{
			\coordinate[shift={(-\dx, 0cm)}] (a\s) at (\anglet + \s: \hairRadius);
			\filldraw (a\s) circle[radius=0.1];
			\draw (lPoint) -- (a\s); 
		}
		\foreach \s in {1, ..., \vertexNumber}
		{
			\pgfmathsetmacro{\angles}{(360/\vertexNumber)* \s + \rotAngleX}
			\coordinate[shift={(-\dx, 0cm)}] (b\s) at (\angles : \kRadius);
			\draw (lPoint) -- (b\s);
		}

		\pgfmathsetmacro{\anglet}{(360 / \vertexNumber)*\t + \rotAngleY}
		\coordinate[shift={(\dx, 0cm)}] (rPoint) at (\anglet : \kRadius);
		\filldraw (rPoint) circle[radius=0.1];
		\foreach \s in {-\hairAngle,\hairAngle}
		{
			\coordinate[shift={(\dx, 0cm)}] (c\s) at (\anglet + \s: \hairRadius);
			\filldraw (c\s) circle[radius=0.1];
			\draw (rPoint) -- (c\s); 
		}
		\foreach \s in {1, ..., \vertexNumber}
		{
			\pgfmathsetmacro{\angles}{(360/\vertexNumber)* \s + \rotAngleY}
			\coordinate[shift={(\dx, 0cm)}] (d\s) at (\angles : \kRadius);
			\draw (rPoint) -- (d\s);
		}

		\draw (lPoint) -- (rPoint);
	}
\end{tikzpicture}
\caption{The graph $(K_{8} \Box K_{2}) \corona \emptyg{2}$ has pretty good state transfer.}
\botline{.15in}
\end{figure}

\begin{thm}
For all $n \ne 1, 3$, the graph $(K_n \Box K_2) \corona \emptyg{2}$ has pretty good state
transfer between the vertices of $\emptyg{2}$.

\begin{proof}
The set of eigenvalues of $K_n \Box K_2$ is given by $\{n, n-2, 0, -2\}$. 
Similar to the proof of Theorem \ref{thm: selflooppgst1}, by Lemma \ref{lem:pgst-necessity}, 
for each $\lambda \in \Sp(K_{n} \Box K_{2})$, we require that
\begin{align}
\cos\left( \frac{t}{2} \lambda \right) \cos\left( \frac{t}{2}\sqrt{\lambda^2 + 8} \right) \approx -1
\end{align}
If $n = 1$ or $3$, this is impossible since $1$ is an eigenvalue.

For all other $n$, we choose times of the form $t = 4\pi\ell$ so that the
first cosine is always $1$ and we simultaneously approximate
\begin{equation} \label{eqn: bunkbed-system}
	\begin{array}{rl}
	2 \ell \sqrt{n^2 + 8} & \approx \ 2r + 1, \\
	2 \ell \sqrt{(n-2)^2 + 8} & \approx \ 2s + 1, \\
	4 \ell \sqrt{2} & \approx \ 2p + 1, \\
	4 \ell \sqrt{3} & \approx \ 2q + 1
	\end{array}
\end{equation}
for some integers $\ell, r, s, p$, and $q$. 
If the set $S = \{\sqrt{n^2 + 8}, \sqrt{(n-2)^2 + 8}, \sqrt{3}, \sqrt{2}, 1\}$ is linearly independent over $\Q$,
then Kronecker's Theorem yields an integer $\ell$ which satisfies \eqref{eqn: bunkbed-system}.

We will now consider when $S$ is linearly dependent. 
First note that a number of the form $\sqrt{m^2 + 8}$ is rational only for $m = \pm 1$. 
Since $n \ne 1,3$, both $\sqrt{n^2 + 8}$ and $\sqrt{(n-2)^2 + 8}$ are irrational.
By Lemma \ref{lem: integersquareroot}, we are interested in only three cases:
\begin{align}
	\sqrt{m^2 + 8} 
	\ = \
	\begin{cases}
		a\sqrt{2} \\
		a\sqrt{3} \\
		\frac{a}{b} \sqrt{(m')^2 + 8}
	\end{cases}
\end{align}
for  integers $a$ and $b$,  and distinct integers $m,m' \in \{n, n-2\}$.

For the cases when $\sqrt{m^2 + 8} = a\sqrt{2}$ or $\sqrt{m^2 + 8} = a\sqrt{3}$,
after squaring and viewing both sides modulo $4$, we see that
$a$ must be twice an odd integer, so \eqref{eqn: bunkbed-system} can still be approximated.

Now consider the case when $\sqrt{n^2 + 8}$ and $\sqrt{(n-2)^2 + 8}$ are linearly dependent. 
Again by Lemma \ref{lem: integersquareroot}, $n^2+8$ and $(n-2)^2+8$ must
have the same square-free part. Denote it by $\Delta$ so that the dependence relation can
be expressed by
\begin{subequations}
\begin{align}
\label{eq: dependence1}
	n^2 + 8 = a^2 \Delta, \\ 
	(n-2)^2 + 8 = b^2 \Delta.
\end{align}
\end{subequations}
Subtracting the two equations:
\begin{align}
	(a^2 - b^2)\Delta = 4(n-1).
\end{align}
Because $\Delta$ is square-free, $a^2 - b^2$ must be even and therefore divisible by 4.
Let $4x = a^2 - b^2$. Then making the substitution $n = x\Delta + 1$ into equation
\eqref{eq: dependence1} gives the equation
\begin{align}
	x^2 \Delta^2 + 2x\Delta + 9 = a^2 \Delta.
\end{align}
It follows that $\Delta$ divides $9$. We know that $\Delta \ne 1$, since $\sqrt{n^2 + 8}$ is irrational,
and since $\Delta$ is square-free, it must be that $\Delta = 3$. This reduces to the case when
$\sqrt{m^2 + 8} = a\sqrt{3}$, which we have already handled. This completes the proof.
\end{proof}
\end{thm}


\subsection{Thorny Graphs}

The corona product of a graph with $K_{1}$ is called a {\em thorny} graph (see Gutman \cite{g98}).

\begin{thm} \label{thm: odd-cube-pgst}
Let $G$ be a graph and let $u,v$ be two of its vertices.
Suppose there is perfect state transfer between $u$ and $v$ at time $t = \pi/g$, for some positive integer $g$,
and that $0$ is not in the eigenvalue support of $u$. Then there is pretty good state transfer 
between $(u,0)$ and $(v,0)$ in $G \corona K_1$.

\begin{proof}
Let $S$ be the eigenvalue support of $u$ in $G$.
By Theorem \ref{thm: coutinho-pst-conditions},
if $G$ has perfect state transfer at time $\pi/g$ between the vertices $u$ and $v$, for some integer $g$, 
all eigenvalues in $S$ must be integers. 
For each eigenvalue $\lambda$ in $S$, let $c_\lambda$ be the square-free part of $\lambda^2 + 4$, 	
so that $\Damp_\lambda = \sqrt{\lambda^2 + 4} = s_\lambda \sqrt{c_\lambda}$ for some integers $s_\lambda$. 
Since $0$ is not in the eigenvalue support of $u$, 
then $\Damp_\lambda$ is irrational and $c_\lambda > 1$ for each $\lambda$ in $S$.
By Lemma \ref{lem: integersquareroot},
\begin{equation}
	\{ \sqrt{c_\lambda} : \lambda \in \supp_G(u)\} \cup \{1\}
\end{equation}
is linearly independent over $\mathbb Q$. 
Kronecker's Theorem implies that there exist integers $\ell, q_\lambda$ such that
\begin{equation}
	\ell \sqrt{c_\lambda} - q_\lambda \ \approx \ -\frac{\sqrt{c_\lambda}}{2g}.
\end{equation}
Multiplying by $4s_\lambda$ yields that
\begin{equation}
	\left(4\ell + \frac{2}{g}\right) \Damp_\lambda \ \approx \ 4 q_\lambda s_\lambda.
\end{equation}
Therefore, at $t = (4\ell + 2/g)\pi$, we have $\cos(\Damp_\lambda t /2) \approx 1$ for each $\lambda$ in $S$. 
By Proposition \ref{prop: adj-transition-corona},
\begin{equation}\begin{split}
\ebra{(u,0)} e^{-\ii tA(G\corona K_1)} \eket{(v,0)} 
	& = \sum_{\lambda \in \Sp(G)} e^{-\ii t\lambda/2} \left(\cos(\Damp_\lambda t/2) 
	- i\frac{\lambda}{\Damp_\lambda} \sin(\Damp_\lambda t/2) \right) \ebra{u} \eigenA{\lambda}{G} \eket{v} \\
	& \approx \sum_{\lambda \in \Sp(G)} e^{- \ii (2\pi) \ell \lambda} e^{-\ii\frac{\pi}{g}\lambda} 
		\ebra{u}\eigenA{\lambda}{G} \eket{v} \\
	& = \ebra{u} e^{-\ii (\pi/g) A(G)} \eket{v},
\end{split}\end{equation}
because all the eigenvalues $\lambda$ in the support of $u$ are integers. 
Since $G$ has perfect state transfer between $u$ and $v$ at time $\pi/g$, this completes the proof.
\end{proof}
\end{thm}

When zero is in the eigenvalue support of $u$, we need a slightly stronger condition to 
get pretty good state transfer in $G\corona K_1$.

\begin{thm} \label{thm: even-cube-pgst}
Let $G$ be a  graph having zero as an eigenvalue. Suppose that $G$ has perfect state transfer at time
$\pi/2$ between vertices $u$ and $v$. Then there is pretty good state transfer between 
$(u,0)$ and $(v,0)$ in the corona $G \corona K_1$.

\begin{proof}
Similar to the proof of Theorem \ref{thm: odd-cube-pgst}, for each $\lambda \in \supp_G(u)$, we can write $\Damp_\lambda = s_\lambda \sqrt{c_\lambda}$ where 
$c_\lambda$ is the square-free part of $\lambda^2+4$ and $s_\lambda$ is an integer.
Note that $c_\lambda = 1$ if and only if $\lambda=0$.
The set 
\begin{equation}
\{ \sqrt{c_\lambda} : \lambda \in \supp_G(u), \lambda\neq 0\} \cup \{1\}
\end{equation}
is linearly independent over $\mathbb Q$.

For $\lambda \neq 0$, Kronecker's Theorem implies that there exist integers $l$ and $q_\lambda$ such that
\begin{equation}
	\ell \sqrt{c_\lambda} - q_\lambda \ \approx \ -\frac{\sqrt{c_\lambda}}{4}+\frac{1}{2s_\lambda}.
\end{equation}
At time $t=(4\ell +1)\pi$, we have $\cos(\Damp_0 t/2)=-1$, 
and $\cos(\Damp_\lambda t/2) \approx -1$ for $\lambda \neq 0$.
By Proposition \ref{prop: adj-transition-corona},
\begin{equation}\begin{split}
\ebra{(u,0)} e^{-\ii tA(G\corona K_1)} \eket{(v,0)} 
	& = \sum_{\lambda \in \Sp(G)} e^{-\ii t\lambda/2} \left(\cos(\Damp_\lambda t/2) 
	- i\frac{\lambda}{\Damp_\lambda} \sin(\Damp_\lambda t/2) \right) \ebra{u} \eigenA{\lambda}{G} \eket{v} \\
	& \approx - \sum_{\lambda \in \Sp(G)} e^{- \ii (2\pi) \ell \lambda} e^{-\ii\frac{\pi}{2}\lambda} 
		\ebra{u}\eigenA{\lambda}{G} \eket{v} \\
	& = - \ebra{u} e^{-\ii (\pi/2) A(G)} \eket{v}.
\end{split}\end{equation}
Since $G$ has perfect state transfer between $u$ and $v$ at time $\pi/2$, this completes the proof.

\end{proof}
\end{thm}

In contrast to Corollary \ref{cor: drg-corona-no-periodic}, the following shows that 
certain thorny distance-regular graphs have pretty good state transfer, but not perfect state transfer. 

\begin{cor}
Let $G$ be a graph from one of the following families:
\begin{itemize}
\item $d$-cubes $Q_d$, for $d \ge 2$.
	\item Cocktail party graphs $\overline{nK_2}$ when $n$ is even.
	\item Halved $2d$-cubes $\frac{1}{2}Q_{2d}$, for $d \ge 1$.
\end{itemize}
Then $G \corona K_{1}$ has pretty good state transfer between antipodal vertices in $G$.
\end{cor}

For a distance-regular graph $G$ with diameter $d$, let $G_{\ell}$ be a graph obtained from $G$
by connecting two vertices $u$ and $v$ if they are at distance $\ell$ from each other, where
$\ell$ ranges from $0$ to $d$. It is customary to denote $A_{\ell}(G)$ as the adjacency matrix of graph $G_{\ell}$.
We say $G$ is {\em antipodal} if $G_{d}$ is a disjoint union of cliques of the same size; here, these cliques
are called the antipodal {\em classes} or {\em fibres} of $G$.

\begin{lem}[Coutinho \etal \cite{cggv15}, Lemma 4.4]  \label{lem: antipodal-class-decomposition}
Let $G$ be a distance-regular graph of diameter $d$ with eigenvalues $\lambda_0 > \lambda_1 > \ldots > \lambda_d$.
Let the spectral decomposition of the adjacency matrix of $G$ be given by
$A(G) = \sum_{j=0}^d \lambda_j \eigenA{j}{G}$.
Suppose that $G$ is antipodal with classes of size two. 
Then
\begin{equation}
	A_d(G) \eigenA{j}{G} = (-1)^j \eigenA{j}{G}.
\end{equation}
\end{lem}

\newcommand*{\qRadius}{1}
\newcommand*{\kscale}{2/3}
\pgfmathsetmacro{\scale}{sin(90 + 45/2)/sin(45/2)}
\begin{figure}[H]
	\centering
	\topline{.15in}
	\begin{tikzpicture}[scale = 0.6]
		\foreach \t in {1,...,8}
		{
		\pgfmathsetmacro{\angle}{45 * \t}
			\filldraw (\angle : \qRadius) circle[radius=0.1];
			\filldraw  (\angle: \scale*\qRadius) circle[radius=0.1];
			\draw (\angle : \scale*\qRadius) -- ({\angle + 45}: \scale*\qRadius);
			\draw (\angle: \qRadius) -- ({\angle + 135}: \qRadius);
			\foreach \s in {-1, 1}
			{
				\draw (\angle: \scale*\qRadius) -- ({\angle + \s*45}: \qRadius);
			}
			\filldraw (\angle: {\qRadius*(1 + \kscale)}) circle[radius=0.1];
				\draw (\angle: \qRadius) -- +(\angle: \qRadius*\kscale);
			\filldraw (\angle: {\qRadius*(\scale + \kscale)}) circle[radius=0.1];
				\draw (\angle: \qRadius*\scale) -- +(\angle: \qRadius*\kscale);
		}		
		\filldraw[color=black, fill=white] (0: \scale*\qRadius) circle[radius=0.1];
		\filldraw[color=black, fill=white] (180: \scale*\qRadius) circle[radius=0.1];
	\end{tikzpicture}
	\caption{$Q_4 \corona K_1$ has PGST between the white vertices.}
	\botline{.15in}
\end{figure}
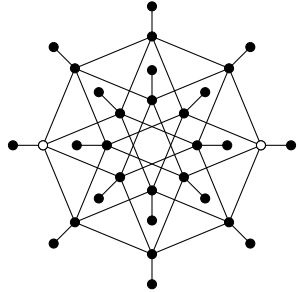

The following observation shows that it is not necessary for $G$ to have perfect state transfer 
in order for $G \corona H$ to have pretty good state transfer.

\begin{thm}
Let $n \ge 3$ be an odd integer, and let $u$ and $v$ be antipodal vertices of the cocktail party graph $\overline{nK_2}$. 
Then there is pretty good state transfer between $(u,0)$ and $(v,0)$ in $\overline{nK_2} \corona K_1$.

\begin{proof}
The eigenvalues of $\overline{nK_2}$ are
$\lambda_0 = 2n-2$, $\lambda_1 = 0$, and $\lambda_2 = -2$.
Since the cocktail-party graph is an antipodal distance-regular graph with fibers of size two, 
from Lemma \ref{lem: antipodal-class-decomposition}, the eigenprojectors satisfy
\begin{equation}
	\ebra{u} \eigenA{j}{\overline{nK_2}} \eket{v} = (-1)^j \ebra{u} \eigenA{j}{\overline{nK_2}} \eket{u}
\end{equation}
for $j = 0,1,2$.
From Proposition \ref{prop: adj-transition-corona}, 
letting $\Damp_{j} = \sqrt{\lambda_{j}^{2} + 4}$,
it suffices to approximate $e^{-\ii t\lambda_j/2} \approx 1$ and
\begin{equation} \label{eq: odd-cocktail-cosines}
	\cos\left(\Damp_j t/ 2\right) \approx (-1)^{j+1}.
\end{equation}
For all integers $\ell$, at time $t = 4\pi\ell$, we have $e^{-\ii t\lambda_j/2} = 1$ and $\cos(\Damp_1 t /2) = 1$.
We will apply Kronecker's Theorem to $\Damp_0$ and $\Damp_2$. Note that $\Damp_2 = 2\sqrt{2}$ while
\begin{equation}
	\Damp_0 = \sqrt{4 + 4(n-1)^2} = 2\sqrt{1 + (n-1)^2},
\end{equation}
and $1 + (n-1)^2$ is always odd when $n$ is odd. 
Let $c_0$ denote the square-free part of $\Damp_0^2$;
then, $\Damp_0 = 2s_0\sqrt{c_0}$ for some odd integer $s_0$. Choose integers $\ell$, $q_0$, and $q_2$ such that
\begin{align}
	\ell \sqrt{c_0} - q_0 \approx \frac{1}{4}, \\
	\ell \sqrt{2} - q_2 \approx \frac{1}{4}.
\end{align}
At time $t = 4\pi \ell$, we note $t \Damp_0/2 \approx 4\pi q_0 s_0 + \pi s_0$ and $t\Damp_2/2 \approx 4\pi q_2+\pi$, 
implying that 
\eqref{eq: odd-cocktail-cosines} is satisfied.
\end{proof}
\end{thm}


\section*{Acknowledgments}

We thank Chris Godsil for helpful remarks on state transfer and Adi Makmal for discussions 
about embedded hypercubes.
The research of E.A., Z.B., J.M., and C.T. was supported by NSF grant DMS-1262737
and NSA grant H98230-15-1-0044.
Part of this work was done while C.T. was visiting Institut Henri Poincar\'{e}.


\bibliographystyle{plain}

\begin{thebibliography}{10}

\bibitem{abcmt-laplacian}
{E.} Ackelsberg, {Z.} Brehm, {A.} Chan, {J.} Mundinger, and {C.} Tamon.
\newblock {L}aplacian state transfer in coronas.
\newblock quant-ph/1508.05458.

\bibitem{bfflott12}
{R.} Bachman, {E.} Fredette, {J.} Fuller, {M.} Landry, {M.} Opperman, {C.}
  Tamon, and {A.} Tollefson.
\newblock Perfect state transfer on quotient graphs.
\newblock {\em Quantum Information and Computation}, 12(3\&4):293--313, 2012.

\bibitem{bps07}
{S.} Barik, {S.} Pati, and {B.} Sarma.
\newblock The spectrum of the corona of two graphs.
\newblock {\em SIAM Journal on Discrete Mathematics}, 21:47--56, 2007.

\bibitem{b03}
{S.} Bose.
\newblock Quantum communication through an unmodulated spin chain.
\newblock {\em Physical Review Letters}, 91(20):207901, 2003.

\bibitem{bcn89}
{A. E.} Brouwer, {A. M.} Cohen, and {A. } Neumaier.
\newblock {\em Distance-regular graphs}.
\newblock Springer-Verlag, 1989.

\bibitem{clms09}
{A.} Casaccino, {S.} Lloyd, {S.} Mancini, and {S.} Severini.
\newblock Quantum state transfer through a qubit network with energy shifts and
  fluctuations.
\newblock {\em International Journal of Quantum Information}, 7(8):1417--1427,
  2009.

\bibitem{ccdfgs03}
{A.} Childs, {R.} Cleve, {E.} Deotto, {E.} Farhi, {S.} Gutmann, and {D.}
  Spielman.
\newblock Exponential algorithmic speedup by a quantum walk.
\newblock In {\em Proceedings of the 35th ACM Symposium on Theory of
  Computing}, pages 59--68, 2003.

\bibitem{cddekl05}
{M.} Christandl, {N.} Datta, {T.} Dorlas, {A.} Ekert, {A.} Kay, and {A.}
  Landahl.
\newblock Perfect transfer of arbitrary states in quantum spin networks.
\newblock {\em Physical Review A}, 71:032312, Mar 2005.

\bibitem{cdel04}
{M.} Christandl, {N.} Datta, {A.} Ekert, and {A.} Landahl.
\newblock Perfect state transfer in quantum spin networks.
\newblock {\em Physical Review Letters}, 92:187902, 2004.

\bibitem{c14}
{G.} Coutinho.
\newblock {\em Quantum State Transfer in Graphs}.
\newblock PhD thesis, University of Waterloo, 2014.

\bibitem{cggv15}
{G.} Coutinho, {C.} Godsil, {K.} Guo, and {F.} Vanhove.
\newblock Perfect state transfer on distance-regular graphs and association
  schemes.
\newblock {\em Linear Algebra and its Applications}, 478:108--130, 2015.

\bibitem{fg13}
{X.} Fan and {C.} Godsil.
\newblock Pretty good state transfer on double stars.
\newblock {\em Linear Algebra and Its Applications}, 438(5):2346--2358, 2013.

\bibitem{fgg08}
{E.} Farhi, {J.} Goldstone, and {S.} Gutmann.
\newblock A quantum algorithm for the hamiltonian {NAND} tree.
\newblock {\em Theory of Computing}, 4(8):169--190, 2008.

\bibitem{fg98}
{E.} Farhi and {S.} Gutmann.
\newblock Quantum computation and decision trees.
\newblock {\em Physical Review A}, 58:915--928, 1998.

\bibitem{fh70}
{R.} Frucht and {F.} Harary.
\newblock On the corona of two graphs.
\newblock {\em Aequationes Mathematicae}, 4(3):322--325, 1970.

\bibitem{gbs08}
{A.} Ghosh, {S.} Boyd, and {A.} Saberi.
\newblock Minimizing effective resistance of a graph.
\newblock {\em SIAM Review}, 50(1):37--66, 2008.

\bibitem{g11}
{C.} Godsil.
\newblock Periodic graphs.
\newblock {\em Electronic Journal of Combinatorics}, 18(1):P23, 2011.

\bibitem{g11-dm}
{C.} Godsil.
\newblock State transfer on graphs.
\newblock {\em Discrete Mathematics}, 312(1):129--147, 2011.

\bibitem{g12}
{C.} Godsil.
\newblock When can perfect state transfer occur?
\newblock {\em Electronic Journal of Linear Algebra}, 23, 2012.

\bibitem{g15}
{C.} Godsil.
\newblock Graph spectra and quantum walks, 2015.

\bibitem{gkss12}
{C.} Godsil, {S.} Kirkland, {S.} Severini, and {J.} Smith.
\newblock Number-theoretic nature of communication in quantum spin systems.
\newblock {\em Physical Review Letters}, 109(5):050502, 2012.

\bibitem{gr01}
{C.} Godsil and {G.} Royle.
\newblock {\em Algebraic Graph Theory}.
\newblock Springer, 2001.

\bibitem{g98}
{I.} Gutman.
\newblock Distance of thorny graphs.
\newblock {\em Publications de l'Institute Math\'{e}matique}, 63(77):31--36,
  1998.

\bibitem{hw00}
{G. H.} Hardy and {E. M.} Wright.
\newblock {\em An Introduction to the Theory of Numbers}.
\newblock Oxford University Press, fifth edition, 2000.

\bibitem{mzmtb14}
{A.} Makmal, {M.} Zhu, {D.} Manzano, {M.} Tiersch, and {H. J.} Briegel.
\newblock Quantum walks on embedded hypercubes.
\newblock {\em Physical Review A}, 90:022314, Aug 2014.

\bibitem{f60}
{D.} Newman and {H.} Flanders.
\newblock Solution to problem 4797.
\newblock {\em American Mathematical Monthly}, 67(2):188--189, 1960.

\bibitem{vz12}
{L.} Vinet and {A.} Zhedanov.
\newblock Almost perfect state transfer in quantum spin chains.
\newblock {\em Physical Review A}, 86:052319, 2012.

\end{thebibliography}


\newpage
\appendix

\section{Auxiliary lemma}

\ignore{
\newcommand{\Real}{\mathfrak{Re}}

\begin{prop} \label{pst}
Let $z_{1},\ldots,z_{m}$ be complex numbers where $|z_{k}| \le 1$ for each $k$.
Let $\vec{p} \in (0,1)^{m}$ be a probability vector.
Suppose that $|\sum_{k=1}^{m} p_{k}z_{k}| = 1$. Then $|z_{k}| = 1$ for each $k$.
Moreover, there is a real $\theta$ so that $z_{k} = e^{i\theta}$.
\begin{proof}
Let $Z = \sum_{k=1}^{m} p_{k}z_{k}$. Since $|Z| = 1$, we have
\begin{eqnarray} \label{eqn:unimodular}
1 & = & \sum_{k} p_{k}^{2}|z_{k}|^{2} + 2\sum_{k < \ell} p_{k}p_{\ell}\Real(z_{k}\overline{z}_{\ell}) \\
	& \le & \sum_{k} p_{k}^{2}|z_{k}|^{2} + 2\sum_{k < \ell} p_{k}p_{\ell}|z_{k}||\overline{z}_{\ell}| \\
	& \le & 1.
\end{eqnarray}
This shows $|z_{k}| = 1$ for each $k$, since otherwise we derive a contradiction.
So, we may assume that $z_{k} = e^{i\alpha_{k}}$ for some real $\alpha_{k}$.
Revisiting \eqref{eqn:unimodular}, we see that
\begin{equation} \label{eqn:uniphase}
1 = \sum_{k} p_{k}^{2} + 2\sum_{k < \ell} p_{k}p_{\ell}\cos(\alpha_{k}-\alpha_{\ell}) \le 1,
\end{equation}
which shows that $\alpha_{k} = \alpha_{\ell}$ for all $k \neq \ell$.
\end{proof}
\end{prop}
}

\begin{lem} \label{lem:pgst-necessity}
Let $z_{1}(t),z_{2}(t),\ldots,z_{m}(t)$ be complex-valued functions satisfying
\begin{itemize}
\item $z_{1}(t) \equiv 1$ or $z_{1}(t) \equiv -1$.
\item for $j=2,\ldots,m$, 
	\begin{equation} \label{eqn:z-defined}
	z_{j}(t) = e^{\ii \alpha_{j}t}(\cos(\theta_{j} t) + \ii\Delta_{j}\sin(\theta_{j} t)),
	\end{equation}
	for some real-valued constants $\alpha_{j},\theta_{j}$ and $\Delta_{j}$ such that $|\Delta_{j}| < 1$. 
\end{itemize}
Let $(p_{1},\ldots,p_{m}) \in (0,1)^{m}$ be a probability distribution. 
Suppose for each $\epsilon > 0$, there is a $\tau > 0$ such that 
\begin{equation} \label{eqn:pgst-convex-combination}
\left|\sum_{j=1}^{m} p_{j}z_{j}(\tau)\right| > 1 - \epsilon.
\end{equation}
There exists a constant $M$ such that for $\epsilon < \min\{p_{1},\ldots,p_{m}\}/2$ and for each $j=2,\ldots,m$
\begin{equation} \label{eqn:pgst-necessary-condition}
|z_{1}(\tau) - \cos(\theta_{j}\tau)\cos(\alpha_{j}\tau)| < M\sqrt{\epsilon}.
\end{equation}

\begin{proof}
We write the functions $z_{j}(t)$ in polar forms as
\begin{equation} \label{eqn:z-redefined}
z_{j}(t) = r_{j}(t)e^{\ii\beta_{j}(t)}
\end{equation}
for some real-valued function $\beta_{j}(t)$
where 
$r_{1}(t) = 1$ and
\begin{equation} \label{eqn:r-defined}
r_{j}(t) = \sqrt{1 + (\Delta_{j}^{2} - 1)\sin^{2}(\theta_{j}t)}, 
\ \hspace{.2in} \
\mbox{ for $j=2,\ldots,m$. }
\end{equation}
Since $|\Delta_{j}| < 1$, note $0 \le r_{j}(t) \le 1$ for $j=1,\ldots,m$.

Using triangle inequality on \eqref{eqn:pgst-convex-combination}, we get
\begin{equation}
1 - \epsilon
	\ < \ \sum_{k=1}^{m} p_{k}r_{k}(\tau) 
	\ \le \ \sum_{k=1}^{m} p_{k} + p_{j}(r_{j}(\tau) - 1)
	\ = \ 1 - p_{j} + p_{j}r_{j}(\tau), \qquad \text{for $j=1,\ldots,m$.}
\end{equation}
Hence $r_{j}(\tau) > (p_{j} - \epsilon)/p_{j}$ for $j=1,\ldots,m$.
For the rest of the proof, we assume $\epsilon < p_{\star}/2$ where
$p_{\star} = \min\{p_{1},\ldots,p_{m}\}$. As a result,
\begin{equation} \label{eqn:r-lowerbound}
r_{j}(\tau) \ > \ 1 - \frac{\epsilon}{p_{j}} \ > \ \frac{1}{2}
\ \hspace{.2in} \
\mbox{ for $j=1,\ldots,m$. }
\end{equation}
Squaring the left-hand side of \eqref{eqn:pgst-convex-combination}, we get from \eqref{eqn:z-redefined} and \eqref{eqn:r-lowerbound} that
\begin{eqnarray}
\left|\sum_{k=1}^{m} p_{k}z_{k}(\tau)\right|^{2}
	& = & 
	\sum_{k=1}^{m} p_{k}^{2}r_{k}(\tau)^{2}
	+ 2\sum_{k < \ell} p_{k}p_{\ell}r_{k}(\tau)r_{\ell}(\tau)\cos(\beta_{k}(\tau)-\beta_{\ell}(\tau)) \\
	& \le &
	\left(\sum_{k=1}^{m} p_{k}\right)^{2}
	+ 2p_{1}p_{j}r_{1}(\tau)r_{j}(\tau)[\cos(\beta_{1}(\tau)-\beta_{j}(\tau)) - 1]\\
	& \le &
	1
	- 2p_{1}p_{j}r_{j}(\tau)|\cos(\beta_{1}(\tau)-\beta_{j}(\tau)) - 1| \\
	& \le &
	\label{eqn:one-minus}
	1
	- p_{1}p_{j}|\cos(\beta_{1}(\tau)-\beta_{j}(\tau)) - 1|.
\end{eqnarray}
Since $\cos\beta_{1}(\tau) = z_{1}(\tau) \in \{-1,+1\}$, we have
\begin{equation} \label{eqn:pm-one}
|\cos(\beta_{1}(\tau) - \beta_{j}(\tau)) - 1|
=
|z_{1}(\tau)\cos\beta_{j}(\tau) - 1|
=
|\cos\beta_{j}(\tau) - z_{1}(\tau)|.
\end{equation}
It follows from \eqref{eqn:pgst-convex-combination}, \eqref{eqn:one-minus} and \eqref{eqn:pm-one} that,
for $j=2,\ldots,m$,
\begin{equation} \label{eqn:m1}
|\cos\beta_{j}(\tau) - z_{1}(\tau)| 
\ < \ 
\frac{\epsilon(2-\epsilon)}{p_{1}p_{j}}.
\end{equation}
From \eqref{eqn:z-defined} and \eqref{eqn:z-redefined}, we get
\begin{eqnarray} \label{eqn:almost-diff}
|\cos\beta_{j}(\tau) - z_{1}(\tau)|
	& \ge &
	\left|\frac{\cos(\alpha_{j}\tau)\cos(\theta_{j}\tau)}{r_{j}(\tau)} - z_{1}(\tau)\right|
		- \left|\frac{\Delta_{j}\sin(\alpha_{j}\tau)\sin(\theta_{j}\tau)}{r_{j}(\tau)}\right|.
\end{eqnarray}
Using \eqref{eqn:r-defined}, we get
\begin{equation}
\left|\sin(\theta_{j}(\tau)\right|=\sqrt{\frac{1-r_j(\tau)^2}{1-\Delta_{j}^2}}.
\end{equation}
Together with \eqref{eqn:r-lowerbound}, we bound the second summand in \eqref{eqn:almost-diff} as follows:
\begin{eqnarray}
\left|\frac{\Delta_{j}\sin(\alpha_{j}\tau)\sin(\theta_{j}\tau)}{r_{j}(\tau)}\right|
	& \le & \frac{\Delta_{j}}{\sqrt{1 - \Delta_{j}^{2}}} \sqrt{\frac{1}{r_{j}(\tau)^{2}} - 1} \\
	& < & \frac{\Delta_{j}}{\sqrt{1 - \Delta_{j}^{2}}} \frac{\sqrt{(2p_{j}-\epsilon)\epsilon}}{(p_{j}-\epsilon)} 
	\ =: \ M'(j).
\end{eqnarray}
As a result, we get
\begin{eqnarray}
|\cos\beta_{j}(\tau) - z_{1}(\tau)|
	& > & 
	\left|z_{1}(\tau) - \frac{1}{r_{j}(\tau)}\cos(\alpha_{j}\tau)\cos(\theta_{j}\tau)\right| - M'(j) \\
	& \ge & 
	|z_{1}(\tau) - \cos(\alpha_{j}\tau)\cos(\theta_{j}\tau)| - 
		\left(\frac{1}{r_{j}(\tau)} - 1\right) - M'(j) \\
	& \ge & 
	\label{eqn:m3}
	|z_{1}(\tau) - \cos(\alpha_{j}\tau)\cos(\theta_{j}\tau)| - 
		\frac{\epsilon}{p_{j} - \epsilon} - M'(j) 
\end{eqnarray}
It follows from \eqref{eqn:m1} and \eqref{eqn:m3} that
\begin{eqnarray}
|z_{1}(\tau) - \cos(\alpha_{j}\tau)\cos(\theta_{j}\tau)|
	& < &
	\frac{\epsilon(2-\epsilon)}{p_{1}p_{j}} + \frac{\epsilon}{p_{j}-\epsilon} + M'(j) \\
	& < &
	\frac{\epsilon(2-\epsilon)}{p_{\star}^{2}} + \frac{2\epsilon}{p_{\star}} 
		+ \frac{\Delta_{j}}{\sqrt{1 - \Delta_{j}^{2}}} \frac{2\sqrt{2\epsilon}}{p_{\star}} 
\end{eqnarray}
since $(p_{j} - \epsilon)/p_{j} \ge 1/2$ by \eqref{eqn:r-lowerbound}.
Let
\begin{equation}
M = \max_{j}\left\{ \frac{2}{p_{\star}^{2}} + \frac{2}{p_{\star}}
		+ \frac{\Delta_{j}}{\sqrt{1 - \Delta_{j}^{2}}} \frac{2\sqrt{2}}{p_{\star}} \right\}.
\end{equation}
Then 
\begin{equation}
|z_{1}(\tau) - \cos(\alpha_{j}\tau)\cos(\theta_{j}\tau)|
\ < \ M\sqrt{\epsilon}.
\end{equation}
\end{proof}
\end{lem}


\end{document}